\documentclass[MSNbibl,number,citesort,seceqn,dvips]{arxbj}
\usepackage{graphicx}

\aid{0}
\volume{19}
\issue{4}
\pubyear{2013}
\firstpage{1150}
\lastpage{1175}
\doi{10.3150/12-BEJSP13} 

\makeatletter

\renewcommand{\dvtx}{:}

\newtheorem{theorem}{Theorem}[section]
\newtheorem{lemma}[theorem]{Lemma}
\newtheorem{proposition}[theorem]{Proposition}

\newremark{remark}[theorem]{Remark}

\newcommand{\mc}[1]{{\mathcal #1}}
\newcommand{\bb}[1]{{\mathbb #1}}

\def\la{\longrightarrow}
\def\t{{\mathcal T}}
\def\u{{\mathcal U}}
\def\wt{\widetilde}
\def\wh{\widehat}
\def\be{{\mathbf e}}
\def\bm{{\mathbf m}}
\def\bq{{\mathbf q}}
\def\ov{\overline}
\def\vep{\varepsilon}
\newcommand{\build}[3]{\,\mathop{#1}\limits_{#2}\limits^{#3}\,}

\makeatother

\begin{document}
\begin{frontmatter}

\title{Quadrangulations with no pendant vertices}
\runtitle{Quadrangulations with no pendant vertices}

\begin{aug}
\author[1]{\fnms{Johel} \snm{Beltran}\thanksref{1}\ead[label=e1]{johel.beltran@pucp.edu.pe}}
\and
\author[2]{\fnms{Jean-Fran\c{c}ois} \snm{Le Gall}\corref{}\thanksref{2}\ead[label=e2]{jean-francois.legall@math.u-psud.fr}}
\runauthor{J. Beltran and J.-F. Le Gall} 
\address[1]{PUCP, Av. Universitaria cdra. 18,
San Miguel, Ap. 1761, Lima 100, Per\'u and IMCA,
Calle los Bi\'ologos 245, Urb. San C\'esar Primera Etapa, Lima 12, Per\'u.
\printead{e1}}
\address[2]{D\'epartement de math\'ematiques, Universit\'e Paris-Sud,
91405 Orsay C\'edex,
France.\\ \printead{e2}}
\end{aug}


%
\begin{abstract}
We prove that the metric space associated with a uniformly distributed
planar quadrangulation with $n$ faces and no pendant vertices converges
modulo a suitable rescaling to the Brownian map. This is a first step
towards the extension of recent convergence results for random planar
maps to the case of graphs satisfying local constraints.
\end{abstract}

%

\end{frontmatter}

\section{Introduction}\label{sec1}

Much recent work has been devoted to studying the convergence of
rescaled planar graphs, viewed as
metric spaces for the graph distance, towards
the universal limiting object called the Brownian map. In the present
article, we establish such a limit theorem
in a particular instance of planar maps satisfying local constraints,
namely quadrangulations with no
pendant vertices, or equivalently with no vertices of degree $1$.

Recall that a planar map is a proper embedding of
a finite connected graph in the
two-dimensional sphere, considered up to orientation-preserving
homeomorphisms of the sphere. Loops and multiple edges are a priori allowed
(however in the case of bipartite graphs that we will consider, there
cannot be any loop).
The faces of the map are
the connected components of the complement of edges, and the degree of
a face counts the number of edges that are incident to it, with the
convention that
if both sides of an edge are incident to the same face, this edge is
counted twice
in the degree of the face (alternatively, the degree of a face may be
defined as the number
of corners to which it is incident). Let $p\geq3$ be an integer.
Special cases
of planar maps are $p$-angulations (triangulations if $p=3$, or
quadrangulations if $p=4$)
where each face has degree $p$. For technical reasons, one often considers
rooted planar maps, meaning that there is a distinguished oriented edge,
whose tail vertex is called the root vertex.
Planar maps have been
studied thoroughly in
combinatorics, and they also arise in other areas of mathematics. Large
random planar graphs are of interest in theoretical
physics, where they serve as models of random geometry~\cite{ADJ}, in
particular in the
theory of two-dimensional quantum gravity.

The recent paper~\cite{LGU} has established a general convergence
theorem for rescaled
random planar maps viewed as metric spaces. Let $p\geq3$ be such that
either $p=3$
or $p$ is even.
For every integer $n\geq1$, let ${\mathbf m}_n$ be a random planar map
that is uniformly distributed
over the set of all rooted $p$-angulations with $n$ faces
(when $p=3$ we need to restrict our attention to even values of $n$ so
that this set
is not empty). We denote the vertex set of ${\bm}_n$
by $V(\bm_n)$. We equip $V(\bm_n)$ with
the graph distance ${\mathrm{d}}^{\bm_n}_{\mathrm{gr}}$, and we view $(V(\bm
_n),{\mathrm{d}}^{\bm_n}_{\mathrm{gr}})$ as
a random variable taking values in the space $\bb K$ of isometry classes
of compact metric spaces. We equip $\bb K$ with the
Gromov--Hausdorff distance $d_{GH}$ (see, e.g.,~\cite{BBI}) and note that
$(\bb K,d_{GH})$ is a Polish space. The main result of~\cite{LGU}
states that there
exists a random compact metric space $(\bm_\infty,D^*)$ called the
Brownian map,
which does not depend on $p$,
and a constant $c_p>0$ depending on $p$, such that
%
\begin{equation}
\label{universal} \bigl(V(\bm_n), c_p n^{-1/4}{\mathrm{d}}^{\bm_n}_{\mathrm{gr}}\bigr) \build{\la} {n\to\infty}
{(\mathrm{d})}
\bigl(\bm_\infty, D^*\bigr)
\end{equation}
where the convergence holds in distribution in the space $(\bb K,d_{GH})$.
A precise description of the Brownian map is given below at the
beginning of Section~\ref{Proof-main}. The constants
$c_p$ are known explicitly (see~\cite{LGU}) and in particular
$c_4=(\frac{9}{8})^{1/4}$. We observe that the case $p=4$
of (\ref{universal}) has been obtained independently by Miermont \cite
{Mi-quad}, and that the case $p=3$
solves a question raised by Schramm~\cite{Sch}. Note that the first
limit theorem involving the Brownian map
was given in the case of quadrangulations by Marckert and
Mokkadem~\cite{MaMo}, but in a weaker form than stated in~(\ref{universal}).

\begin{figure}

\includegraphics{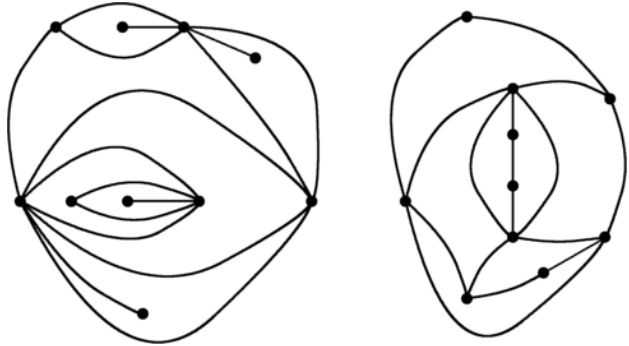}

\caption{Two quadrangulations with $8$ faces. The one on the right is
nice.}\label{fig1}
\end{figure}

In this work, we are interested in planar maps that satisfy additional
local regularity properties.
Under such constraints, one may ask whether the scaling limit is still
the Brownian map, and, if it is,
one expects to get different scaling constants $c_p$.
Note that the general strategy for proving limiting results such as
(\ref{universal}) involves coding
the planar maps by certain labeled trees and deriving asymptotics for
these trees. If the map
is subject to local constraints, say concerning the degree of vertices,
or the absence of multiple edges or of
loops (in the case of triangulations), this leads to certain
conditionings of the trees, which often make
the desired asymptotics much harder to handle. In the present work, we
consider quadrangulations
with no pendant vertices, or equivalently with no vertices of degree
$1$, which we call nice
quadrangulations (see Figure~\ref{fig1}). We let $\mc Q_n^{\mathrm{nice}}$ be the set of
all rooted nice quadrangulations with $n$ faces.
This set is nonempty for every $n\geq2$.

\begin{theorem}
\label{MainTh}
For every $n\geq2$, let $\bq_n$ be uniformly distributed over the set
$\mc Q_n^{\mathrm{nice}}$. Let $V(\bq_n)$ be the vertex set of $\bq_n$ and
let ${\mathrm{d}}^{\bq_n}_{\mathrm{gr}}$
be the graph distance on $V(\bq_n)$. Then,
\[
\bigl(V(\bq_n), \bigl(\tfrac{3}{4}\bigr)^{3/8}n^{-1/4}
{\mathrm{d}}^{\bq_n}_{\mathrm{gr}} \bigr) \build{\la} {n\to\infty} {(\mathrm{d})}
\bigl(\bm_\infty,D^*\bigr)
\]
where $(\bm_\infty, D^*)$ is the Brownian map, and the convergence
holds in distribution in the
space $(\bb K,d_{GH})$.
\end{theorem}

We observe that the limiting space is again the Brownian map, and so
one may say that nice
quadrangulations have asymptotically the same ``shape'' as ordinary
quadrangulations. On the other
hand, the scaling constant is different: Since $(\frac
{3}{4})^{3/8}<(\frac{9}{8})^{1/4}$, distances are
typically larger in nice quadrangulations, as one might have expected.

In relation with Theorem~\ref{MainTh},
we mention the recent work of Bouttier and Guitter~\cite{BG}, which
obtains detailed information
about distances in large quadrangulations with no multiple edges. Note
that a quadrangulation with no multiple edges
is always nice in our sense, but the converse is not true (see the nice
quadrangulation on the right side of Figure~\ref{fig1}).

We view Theorem~\ref{MainTh} as a first step towards the derivation of
similar results in
more difficult cases. A particularly interesting problem is to derive
the analog of (\ref{universal})
for triangulations without loops or multiple edges (type III
triangulations in the
terminology of~\cite{ADJ}). It is known that such a triangulation can be
represented as the tangency graph of a circle packing of the sphere,
and that this representation
is unique up to the conformal transformations of the sphere (the M\"obius transformations).
So assuming that the analog of (\ref{universal}) holds for type III
triangulations, one might expect
to be able to pass to the limit $n\to\infty$ in the associated circle
packings, and to get a canonical embedding
of the Brownian map in the sphere that would satisfy remarkable
conformal invariance properties. One also
conjectures that this canonical embedding would be related to the
recent approach to two-dimensional
quantum gravity which has been developed by Duplantier and Sheffield
\cite{DS} via the Gaussian free field.
The previous questions are among the most fascinating open problems in
the area.

As a final remark, our proofs rely on Schaeffer's bijection between
rooted quadrangulations
and well-labeled trees. One may be tempted to use the version of this
bijection for rooted
and pointed quadrangulations, which avoids the positivity condition on
labels (see, e.g.,
\cite{LGM}). However, the use of this other version of the bijection
in our setting would lead to
certain conditionings (involving the event that the minimal label on
the tree is
attained at two different corners), which seem difficult to handle.

The paper is organized as follows. In Section~\ref{sec2}, we recall Schaeffer's
bijection (we refer to~\cite{CS} for more details) and we identify
those trees that
correspond to nice triangulations. We then state the key limit theorem
for the coding functions of
the random tree associated with a uniformly distributed nice
quadrangulation with $n$ faces. This limit
theorem is the main ingredient of our proof of Theorem~\ref{MainTh} in
Section~\ref{sec3}, which also uses
some ideas introduced in~\cite{LGU} to deal with triangulations. The
proof of the limit theorem for
coding functions is given in Section~\ref{sec4}, which is the most technical
part of the paper.

\section{Trees and quadrangulations}\label{sec2}

\subsection{Labeled trees}\label{sec2.1}

We set $\bb N=\{1,2,\ldots\}$ and by convention $\bb N^0=\{\varnothing
\}$.
We introduce the set
\[
\bb V=\bigcup_{n=0}^\infty\bb
N^n.
\]
An
element
of $\bb V$ is thus a sequence
$u=(u^1,\ldots,u^n)$ of elements of $\bb N$, and we set $|u|=n$, so that
$|u|$ represents the ``generation'' of $u$. If
$u=(u^1,\ldots, u^m)$ and
$v=(v^1,\ldots, v^n)$ belong to $\mathcal{U}$, we write
$uv=(u^1,\ldots, u^m,v^1,\ldots,v^n)$
for the concatenation of $u$ and $v$. In particular,
$u\varnothing=\varnothing u=u$.

If $w\in\bb V$, we write $\bb V_{(w)}$ for the set
of all elements $u\in\bb V$ of the form $u=wv$ for some $v\in\bb V$.
We then set $\bb V^{(w)}= (\bb V \setminus\bb V_{(w)})\cup\{w\}$.

The mapping $\pi\dvtx \bb V\setminus\{\varnothing\}\la\bb V$
is defined by $\pi((u^1,\ldots,u^n))=(u^1,\ldots,u^{n-1})$
($\pi(u)$ is the ``parent'' of $u$).

A plane tree $\tau$ is a finite subset of
$\bb V$ such that:
\begin{enumerate}[(iii)]
\item[(i)] $\varnothing\in\tau$.

\item[(ii)] $u\in\tau\setminus\{\varnothing\}\Rightarrow
\pi(u)\in\tau$.

\item[(iii)] For every $u\in\tau$, there exists an integer $k_u(\tau
)\geq0$
such that, for every $j\in\bb N$, $uj\in\tau$ if and only if $1\leq
j\leq
k_u(\tau)$.
\end{enumerate}

Edges of $\tau$ are all pairs $(u,v)$ where $v\in\tau\setminus\{
\varnothing\}$ and $u=\pi(v)$.
We write $E(\tau)$ for the set of all edges of $\tau$. Every $e\in
E(\tau)$ can therefore
be written as $e=(e_-,e_+)$ where $e_-=\pi(e_+)$.
By definition, the size $|\tau|$ of $\tau$
is the number of edges of $\tau$, $|\tau|=\# E(\tau)=\#\tau-1$.

In what
follows, we see each vertex of the tree $\tau$ as an individual of a
population whose family tree is the tree $\tau$. In (iii) above, the
individuals
of the form $uj$, with $j\in\bb N$, are interpreted as the
``children'' of $u$, and they are ordered
in the obvious way. The number $k_u(\tau)$ is the number of children
of $u$
in $\tau$.
The notions of
an ancestor and a descendant of a vertex $u$ are defined similarly.

\begin{figure}

\includegraphics{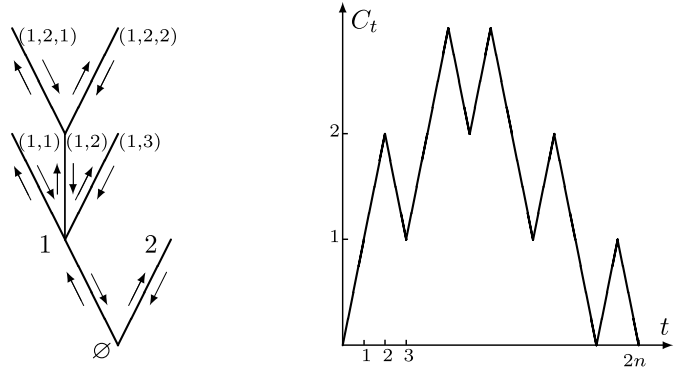}

\caption{A plane tree with $n=7$ edges and its contour function.}
\label{fig2}
\end{figure}

Let $\tau$ be a plane tree and $n=|\tau|$. The contour exploration
sequence of $\tau$ is the finite sequence $v_0,v_1,\ldots,v_{2n}$ which
is defined inductively as follows. First $v_0=\varnothing$, and then,
for every $i\in\{0,\ldots,2n-1\} $, $v_{i+1}$ is either the first child
of $v_{i}$ that does not appear among $v_0,v_1,\ldots,v_{i}$, or, if
there is no such child, the parent of $v_{i}$. Informally, if the tree
is embedded in the plane as suggested in Figure~\ref{fig2}, we imagine the motion
of a particle that starts from the root and traverses the tree from the
left to the right, in the way explained by the arrows of Figure~\ref{fig2}, until
all edges have been explored and the particle has come back to the
root. Then $v_0,v_1,\ldots,v_{2n}$ are the successive vertices visited
by the particle. The contour function of the tree is defined by
$C_i=|v_i|$ for every $i\in\{0,1,\ldots,2n\}$. We extend the function
$C_t$ to the real interval $[0,2n]$ by linear interpolation, and by
convention we set $C_t=0$ for $t\geq2n$. Clearly the tree $\tau$ is
determined by its contour function $(C_t)_{t\geq0}$.

A labeled tree is a pair $(\tau, (U(v))_{v\in\tau})$ that consists
of a plane tree $\tau$ and a collection
$(U(v))_{v\in\tau}$ of integer labels assigned to the vertices of
$\tau$ -- in our
formalism for plane trees, the tree $\tau$ coincides with the set of
all its vertices.
We assume that labels satisfy the following three properties:
\begin{enumerate}[(iii)]
\item[(i)] for every $v\in\tau$, $U(v)\in\bb Z $;
\item[(ii)] $U(\varnothing)=0 $;
\item[(iii)] for every $v\in\tau\setminus\{\varnothing\}$, $U(v)
- U({\pi(v)})\in\{-1,0,1\}$,
\end{enumerate}
where we recall that $\pi(v)$ denotes the parent of $v$. Condition
(iii) just means that when
crossing an edge
of $\tau$ the label can change by at most $1$ in absolute value. We
write $\mc W$
for the set of all labeled trees.

Let $(\tau, (U(v))_{v\in\tau})$ be a labeled tree with $n$ edges. As
we have just seen,
the plane tree $\tau$ is coded by its contour function $(C_t)_{t\geq
0}$. We can similarly
encode the labels by another function $(V_t)_{t\geq0}$, which is
defined as follows.
As above, let $v_0, v_1,v_2,\ldots, v_{2n}$ be the contour exploration
sequence of $\tau$.
We set
\[
V_i=U({v_i})\quad\mbox{for every }i=0,1,\ldots,2n.
\]
Notice that $V_0=V_{2n}=0$.
We extend the function $V_t$ to the real interval $[0,2n]$ by
linear interpolation, and we set $V_t=0$ for $t\geq2n$. We will call
$(V_t)_{t\geq0}$ the ``label function'' of the labeled tree $(\tau,
(U(v))_{v\in\tau})$.
Clearly
$(\tau, (U(v))_{v\in\tau})$ is determined by the pair
$(C_t,V_t)_{t\geq0}$.

We write $\mc W^+$ for the set of all labeled trees with nonnegative
labels (these are sometimes
called well-labeled trees), and for every
$n\geq0$, we write $\mc W^+_n$ for the set of all labeled trees with
$n$ edges in $\mc W^+$.

\subsection{Schaeffer's bijection}\label{sec2.2}

In this section, we fix $n\geq1$ and we briefly recall Schaeffer's
bijection between the set $\mc Q_n$ of
all rooted quadrangulations with $n$ faces and the set $\mc W^+_n$.
We refer to~\cite{CS} for more details. We then characterize
those labeled trees that correspond to nice quadrangulations in this bijection.

To describe Schaeffer's bijection, start from a labeled tree $(\tau,
(U(v))_{v\in\tau})\in\mc W^+_n$, and as above write $v_0,
v_1,v_2,\ldots, v_{2n}$ for the contour exploration sequence of the
vertices of~$\tau$. Notice that each index $i\in\{0,1,\ldots,2n-1\}$
corresponds to exactly one corner of the vertex $v_i$ (a corner of a
vertex $v$ of $\tau$ is an angular sector between two successive edges
of $\tau$ around the vertex $v$). This corner will be called the corner
$i$ in the tree $\tau$.

We extend the
contour exploration sequence periodically, in such
a way that $v_{i+2n}=v_i$ for every integer $i\geq0$. Then, for every
$i\in\{0,1,\ldots,2n-1\}$, we define the successor of $i$ by setting
\[
{\mathrm{succ}}(i)=\cases{\min\bigl\{j\geq i\dvtx  U(v_j)=U(v_i)-1
\bigr\}, & if $U(v_i)>0$,
\cr
\infty, & otherwise.}
\]

\begin{figure}[b]

\includegraphics{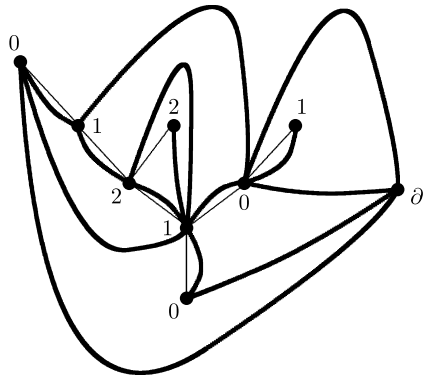}

\caption{Illustration of Schaeffer's bijection. The thin lines
represent the edges of the tree, and the numbers $0,1,\ldots$
are the labels assigned to the different vertices. The thick curves
represent the edges of the associated quadrangulation. The two pendant vertices
are the leaves $v$ such that $U(v)\geq U(\pi(v))$.}
\label{fig3}
\end{figure}

To construct the edges of the quadrangulation associated with $(\tau,
(U(v))_{v\in\tau})$, we proceed
in the following way. We suppose that the tree $\tau$ is drawn in the
plane in the
way suggested in Figure~\ref{fig2}, and we add an extra vertex $\partial$ (outside
the tree).
Then, for every
$i\in\{0,1,\ldots,2n-1\}$,
\begin{enumerate}
\item[$\bullet$] Either $U(v_i)=0$, and we draw an edge between $v_i$
and $\partial$, that starts
from the corner $i$.
\item[$\bullet$] Or $U(v_i)>0$, and we draw an edge between $v_i$ and
$v_{{\mathrm{succ}}(i)}$, that starts
from the corner $i$ and ends at the corner ${\mathrm{succ}}(i)$.
\end{enumerate}
The construction can be made in such a way that the edges do not
intersect, and do not intersect the edges of the tree
(see Figure~\ref{fig3} for an example). The resulting graph,
whose vertex set consists of all vertices of $\tau$ and the vertex
$\partial$, is a quadrangulation with
$n$ faces.
It is rooted at the edge drawn between the vertex $\varnothing$ of
$\tau$ and the vertex $\partial$,
which is oriented in such a way that $\partial$ is the root vertex. We
have thus obtained a
rooted quadrangulation with $n$ faces, which is denoted by $\bq=\Phi
_n(\tau,(U(v))_{v\in\tau})$.
The mapping $\Phi_n$ is Schaeffer's bijection from $\mc W^+_n$ onto
$\mc Q_n$. A key property
of this bijection is the fact that labels on the tree $\tau$ become
distances from the root vertex
$\partial$ in the quadrangulation: If ${\mathrm{d}}^\bq_{\mathrm{gr}}$ stands
for the graph distance on the vertex
set of $\bq$, we have
\[
{\mathrm{d}}^\bq_{\mathrm{gr}}(\partial,v)= U(v) + 1,
\]
for every vertex $v$ of $\tau$ or equivalently for every vertex $v$ of
$\bq$ other than the root vertex.

A leaf of the tree $\tau$ is a vertex with degree $1$. If $v\in\tau
\setminus\{\varnothing\}$,
$v$ is a leaf if and only if $k_v(\tau)=0$, and $\varnothing$ is a
leaf if and only
if $k_\varnothing(\tau)=1$.

\begin{proposition}
\label{nicetree} Let $(\tau, (U(v))_{v\in\tau})\in\mc W^+_n$, and let
$v_0,v_1,\ldots,v_{2n}$ be the contour exploration sequence of $\tau$.
Then the quadrangulation $\Phi_n(\tau,(U(v))_{v\in\tau})$ is nice if
and only if the following two conditions hold.
\begin{enumerate}[(ii)]
\item[(i)] For every leaf $v$ of $\tau$, if $w$ is the (unique)
vertex adjacent to $v$
in the tree $\tau$, we have $U(v)=U(w)-1$.
\item[(ii)] There exists at least one index $i\in\{1,\ldots,2n-1\}$
such that $U(v_i)=0$.
\end{enumerate}
\end{proposition}

Notice that we have always $U(v_0)=U(\varnothing)=0$. Condition (ii)
can be restated by saying that
there are at least two corners of the tree $\tau$ with label $0$. In
particular this condition holds
if $k_\varnothing(\tau)\geq2$.

\begin{pf*}{Proof of Proposition~\ref{nicetree}}
Let us explain why conditions (i) and (ii) are necessary. If (ii) does
not hold, there is only one edge
incident to $\partial$. If there exists a leaf $v$ for which the
property stated in (i) fails, then the only edge
incident to $v$ will be the edge connecting the unique corner of $v$ to
its successor. Conversely,
it is also very easy to check that if conditions (i) and (ii) hold then
every vertex of $\tau$
will be incident to at least $2$ edges in the quadrangulation $\Phi
_n(\tau,(U(v))_{v\in\tau})$: In particular
if $v$ is a leaf of $\tau$ and if $w$ is the vertex adjacent to $v$,
then the successor of one of the corners
of $w$ will be the (unique) corner of $v$. We leave the details to the reader.
\end{pf*}

\begin{remark}
For a general quadrangulation, each leaf $v\ne\varnothing$ of the
associated labeled tree such
that $U(v)\geq U(\pi(v))$ corresponds to a pendant vertex (see Figure~\ref{fig3}).
Using this observation, it
is not hard to prove that
a quadrangulation with $n$ faces has typically about $n/3$ pendant vertices.
\end{remark}

We write $\mc W^{\mathrm{nice}}_n$ for the set of all labeled trees in $\mc
W^+_n$ that satisfy
both conditions in Proposition~\ref{nicetree}.

\subsection{Scaling limits for coding functions}\label{sec2.3}
\label{scalicofu}

In this section, we state the key theorem giving scaling limits for the contour
and label functions of the labeled tree associated with a uniformly distributed
nice quadrangulation with $n$ faces. We first need to introduce the limiting
processes that will appear in this theorem.

We let $\be=(\be_t)_{t\in[0,1]}$ denote a normalized Brownian excursion.
The process $\be$ is defined on a probability space $(\Omega,\mathcal
{A},\bb P)$.
We
consider another real-valued process $Z=(Z_t)_{t\in[0,1]}$ defined
on the same probability space and
such that, conditionally on $\be$, $Z$ is a centered Gaussian process
with covariance
\[
{\bb E}[Z_sZ_t \mid \be]= \min_{r\in[s\wedge t,s\vee t]}
\be_r.
\]
We may and will assume that $Z$ has continuous sample paths. The process
$Z$ can be interpreted as the head of the standard Brownian snake driven
by $\be$.

It is not hard to verify that the distribution of
\[
\min_{t\in[0,1]} Z_t
\]
has no atoms, and that the topological support of this distribution is
$(-\infty,0]$.
Consequently, we can consider for every $r>0$ a process $(\be^{(r)},Z^{(r)})$
whose distribution is the conditional distribution of $(\be,Z)$
knowing that
\[
\min_{t\in[0,1]} Z_t>-r,
\]
and the distribution of $(\be^{(r)},Z^{(r)})$ depends continuously on
$r>0$. Here the
distribution of $(\be^{(r)},Z^{(r)})$ is a probability measure on the
space $C([0,1],\bb R^2)$
of all continuous functions from $[0,1]$ into $\bb R^2$, and
``continuously'' refers to the usual
weak convergence of probability measures.
It is proved in~\cite{LGW}, Theorem 1.1, that we can define a process
$(\be^{(0)},Z^{(0)})$
such that
\[
\bigl(\be^{(r)},Z^{(r)}\bigr)\build{\la} {r\to0} {(\mathrm{d})}
\bigl(\be^{(0)},Z^{(0)}\bigr)
\]
where the convergence holds in distribution in the space $C([0,1],\bb R^2)$.

The following theorem is the key ingredient of the proof of our main result.

\begin{theorem}
\label{conv-coding}
Let $C^{(n)}$ and $V^{(n)}$ be, respectively, the contour function and
the label
function of a random labeled tree distributed uniformly over $\mc
W^{\mathrm{nice}}_n$.
Then,
\[
\bigl(12^{-1/4} n^{-1/2} C^{(n)}_{2nt},
\bigl(\tfrac{3}{4}\bigr)^{3/8} n^{-1/4}
V^{(n)}_{2nt} \bigr)_{0\leq t\leq1} \build{\la} {n\to\infty} {
(\mathrm{d})} \bigl(\be^{(0)}_t,Z^{(0)}_t
\bigr)_{0\leq
t\leq1},
\]
where the convergence holds in distribution in $C([0,1],\bb R^2)$.\vadjust{\goodbreak}
\end{theorem}

The proof of Theorem~\ref{conv-coding} is given in Section~\ref{sec4} below.

\section{Proof of the main theorem}\label{sec3}
\label{Proof-main}

In this section, we explain how to derive Theorem~\ref{MainTh} from
the convergence
of coding functions stated in Theorem~\ref{conv-coding}. Much of what
follows is similar to the arguments
of~\cite{IM}, Section 3, or of~\cite{LGM}, Section 6.2, but we will
provide some details for the sake of completeness.

We start by recalling the definition of the Brownian map.
The first ingredient is the Continuum Random Tree or CRT, which is
conveniently defined as the
tree coded by the Brownian excursion (Aldous~\cite{Al3}). Recall that
if $g\dvtx [0,1]\la\bb R_+$ is a continuous function such that
$g(0)=g(1)=0$, one introduces the equivalence relation on $[0,1]$
defined by
\[
s\sim_g t \quad\mbox{if and only if}\quad g(s)=g(t)=
m_g(s,t),
\]
where $m_g(s,t)=\min\{g(r)\dvtx s\wedge t\leq r\leq s\vee t \}$, and the
tree coded by $g$ is the quotient
space $\t_g:=[0,1] / \sim_g$, which is equipped with the distance
induced by the pseudo-metric
\[
\delta_g(s,t)= g(s)+g(t)-2 m_g(s,t), \quad s,t
\in[0,1].
\]
We write $p_g\dvtx [0,1] \la\t_g$ for the canonical projection. By
convention, $\t_g$ is rooted at $p_g(0)=p_g(1)$.
The CRT $\t_\be$ is then the (random) tree coded by the normalized
Brownian excursion $\be$.

From the definition of the process $Z$, one easily checks that $\bb
E[(Z_s-Z_t)^2 \mid \be]=\delta_\be(s,t)$,
and it follows that we have $Z_s=Z_t$ for every $s,t\in[0,1]$ such
that $s\sim_\be t$, a.s. Hence we may
and sometimes will view $Z$ as indexed by $\t_\be$ rather than by
$[0,1]$. For $a\in\t_\be$,
we interpret $Z_a$ as the ``label'' of the vertex $a$.

We now explain how a trajectorial transformation of $(\be,Z)$ yields
a pair $(\ov\be, \ov Z)$ having the same distribution as $(\be
^{(0)},Z^{(0)})$. By~\cite{LGW}, Proposition 2.5
(and an obvious scaling argument) there exists an a.s. unique time
$s_*\in[0,1]$ such that
$Z_{s_*}=\min\{Z_s\dvtx s\in[0,1]\}$. We then set, for every $t\in[0,1]$,
\begin{eqnarray*}
\ov\be_t &=& \delta_\be(s_*,s_*\oplus t)=
\be_{s_*} + \be_{s_*\oplus t} - 2 m_\be(s_*,s_*\oplus t),
\\
\ov Z_t &=& Z_{s_*\oplus t} -Z_{s_*},
\end{eqnarray*}
where $s_*\oplus t=s+t$ if $s_*+t\leq1$ and $s_*\oplus t=s+t-1$
otherwise. By~\cite{LGW}, Theorem 1.2, the
pair $(\ov\be, \ov Z)$ has the same distribution as $(\be^{(0)},Z^{(0)})$.

One easily verifies that the
property $s_*\oplus t\sim_{\be} s_*\oplus t'$ holds if and only if
$t\sim_{\ov\be} t'$, for
every $t,t'\in[0,1]$, a.s., and it follows that we have $\ov Z_t=\ov
Z_{t'}$ if $t\sim_{\ov\be} t'$.
Hence, we may again view $\ov Z$ as indexed by the tree $\t_{\ov\be}$.

The mapping $t\la s_*\oplus t$ induces an isometry $\mathcal{I}$ from
$\t_{\ov\be}$ onto $\t_\be$, that maps
the root of $\t_{\ov\be}$ to the vertex $p_\be(s_*)$ with minimal
label in $\t_\be$. Furthermore,
we have $\ov Z_a= Z_{\mathcal I(a)}- \min Z$ for every $a\in\t_{\ov
\be}$. To summarize the preceding
discussion, $\t_{\ov\be}$ can be viewed as $\t_\be$ ``re-rooted''
at the vertex with minimal label, and the
labels $\ov Z$ on $\t_{\ov\be}$ are derived from the labels $Z$ on
$\t_\be$ by subtracting the minimal label.

Next, for every $s,t\in[0,1]$, we set
\[
D^\circ(s,t):=\ov Z_s + \ov Z_t -2 \min
_{s\wedge t\leq r\leq s\vee t} \ov Z_r
\]
and, for every $a,b\in\t_{\ov\be}$,
\[
D^\circ(a,b):= \inf\bigl\{ D^\circ(s,t)\dvtx  s,t\in[0,1],
p_{\ov\be
}(s)=a, p_{\ov\be}(t)=b\bigr\}.
\]
Finally, we define a pseudo-metric $D^*$ on $\t_{\ov\be}$ by setting
\[
D^*(a,b):=\inf\Biggl\{\sum_{i=1}^k
D^\circ(a_{i-1},a_i) \Biggr\}
\]
where the infimum is over all choices of the integer $k\geq1$ and of
the finite sequence
$a_0,a_1,\ldots,a_k$ such that $a_0=a$ and $a_k=b$. We set $a\approx
b$ if
and only if $D^*(a,b)=0$ (according to~\cite{IM}, Theorem~3.4, this
holds if and only
if $D^\circ(a,b)=0$).

The Brownian map is the quotient space $\bm_\infty:=\t_{\ov\be}
/ \approx$, which is equipped
with the distance induced by $D^*$. The reader may have noticed that
our presentation is consistent
with~\cite{IM}, but slightly differs from the introduction of \cite
{LGU}, where the Brownian map is
constructed directly from the pair $(\be,Z)$, rather than from $(\ov
\be,\ov Z)$. The previous discussion
about the relations between the trees $\t_\be$ and $\t_{\ov\be}$,
and the labels on these trees, however
shows that both presentations are equivalent. In the present work,
because our limit theorem for the
coding functions of discrete objects involves a pair distributed as
$(\ov\be,\ov Z)$, it will be more
convenient to use the presentation above.

Let us turn to the proof of Theorem~\ref{MainTh}. We let $(\tau
_n,(U_n(v))_{v\in\bq_n})$ be the labeled tree
associated with $\bq_n$, which is uniformly distributed over $\mc
W_n^{\mathrm{nice}}$. As in Theorem
\ref{conv-coding}, we denote the contour function and the label
function of $(\tau_n,(U_n(v))_{v\in\bq_n})$
by $C^{(n)}$ and $V^{(n)}$, respectively. We also write
$u^n_0,u^n_1,\ldots,u^n_{2n}$ for the
contour exploration sequence of $\tau_n$. We then set, for every
$i,j\in\{0,1,\ldots,2n\}$,
\[
{\mathrm{d}}_n(i,j) = {\mathrm{d}}^{\bq_n}_{\mathrm{gr}}
\bigl(u^n_i,u^n_j\bigr)
\]
where ${\mathrm{d}}^{\bq_n}_{\mathrm{gr}}$ stands for the graph distance on
$V(\bq_n)$ (here and in what follows, we use
Schaeffer's bijection to view the vertices of $\tau_n$ as vertices of
$\bq_n$). We extend the definition
of ${\mathrm{d}}_n(i,j)$ to noninteger values of $i$ and $j$ by setting, for
every $s,t\in[0,2n]$,
\begin{eqnarray*}
{\mathrm{d}}_n(s,t)&=&\bigl(s-\lfloor s\rfloor\bigr) \bigl(t-\lfloor t
\rfloor\bigr){\mathrm{d}}_n\bigl(\lceil s\rceil,\lceil t\rceil\bigr) +
\bigl(s-\lfloor s\rfloor\bigr) \bigl(\lceil t\rceil-t\bigr) {\mathrm{d}}_n
\bigl(\lceil s\rceil,\lfloor t\rfloor\bigr)
\\
&&{}+ \bigl(\lceil s\rceil-s\bigr) \bigl(t-\lfloor t\rfloor\bigr){\mathrm{d}}_n\bigl(\lfloor s\rfloor,\lceil t\rceil\bigr) + \bigl(\lceil s
\rceil-s\bigr) \bigl(\lceil t\rceil-t\bigr){\mathrm{d}}_n\bigl(\lfloor s
\rfloor,\lfloor t\rfloor\bigr),
\end{eqnarray*}
where $\lceil t\rceil:=\min\{k\in\bb Z\dvtx  k>t\}$. The same arguments
as in~\cite{IM}, Proposition 3.2,
relying on the bound
%
\begin{equation}
\label{bd-dist} {\mathrm{d}}_n(i,j)\leq V^{(n)}_i +
V^{(n)}_i - 2\min_{i\wedge i\leq k\leq
i\vee j}
V^{(n)}_k +2
\end{equation}
(see~\cite{IM}, Lemma 3.1) and on Theorem~\ref{conv-coding} show that
the sequence
of the laws of the processes $(n^{-1/4}{\mathrm{d}}_n(2ns,2nt))_{0\leq
s,t\leq1}$ is tight
in the space of all probability measures on $C([0,1]^2,\bb R_+)$. Using
this tightness property
and Theorem~\ref{conv-coding}, we can find a
sequence of integers $(n_k)_{k\geq1}$ converging to $+\infty$ and a
continuous random process $(D(s,t))_{0\leq s\leq t}$ such that, along
the sequence $(n_k)_{k\geq1}$,
we have the joint convergence in distribution in $C([0,1]^2,\bb R^3)$,
%
\begin{eqnarray}
\label{cd-dist} && \bigl( 12^{-1/4} n^{-1/2}C^{(n)}_{2nt},
\bigl(\tfrac{3}{4}\bigr)^{3/8} n^{-1/4}V^{(n)}_{2nt},
\bigl(\tfrac{3}{4}\bigr)^{3/8}n^{-1/4} {\mathrm{d}}_n(2ns,2nt) \bigr)_{0\leq s, t\leq1}
\\
&&\quad\build{\la} {n\to\infty} {(\mathrm{d})} \bigl(\ov\be_t,\ov
Z_t, D(s,t) \bigr)_{0\leq s,t\leq1}.\nonumber
\end{eqnarray}
By Skorokhod's representation theorem, we may assume that this
convergence holds a.s.
Passing to the limit $n\to\infty$ in (\ref{bd-dist}), we get that
$D(s,t)\leq D^\circ(s,t)$
for every $s,t\in[0,1]$, a.s. Also, from the fact that ${\mathrm{d}}^{\bq
_n}_{\mathrm{gr}}(\partial, u^n_i)= U(u^n_i) + 1 = V^{(n)}_i +1$
we immediately obtain that $D(0,t)=\ov Z_t$ for every $t\in[0,1]$, a.s.

Clearly, the function $(s,t)\la D(s,t)$ is symmetric, and it also
satisfies the triangle inequality
because ${\mathrm{d}}_n$ does. Furthermore, the fact that ${\mathrm{d}}_n(i,j)=0$
if $u^n_i=u^n_j$ easily implies
that $D(s,t) =0$ for every $s,t\in[0,1]$ such that $s\sim_{\ov\be}
t$, a.s. (see the proof of
Proposition 3.3(iii) in~\cite{IM}). Hence, we may view $D$ as a
random pseudo-metric on $\t_{\ov\be}$.
Since $D\leq D^\circ$ and $D$ satisfies the triangle inequality, the
definition of
$D^*$ immediately shows that $D(a,b)\leq D^*(a,b)$ for every $a,b\in\t
_{\ov\be}$, a.s.

\begin{lemma}
\label{crucial-main}
We have $D(a,b)=D^*(a,b)$ for every $a,b\in\t_{\ov\be}$, a.s.
\end{lemma}

We postpone the proof of Lemma~\ref{crucial-main}
to the end of the section and complete the proof of Theorem~\ref{MainTh}.
We define a correspondence between the metric spaces $(V(\bq
_n)\setminus\{\partial\}, (\frac{3}{4})^{3/8}n^{-1/4}
{\mathrm{d}}^{\bq_n}_{\mathrm{gr}})$ and $(\bm_\infty, D^*)$ by setting
\[
{\mc R}_n:= \bigl\{ \bigl(u^n_{\lfloor2nt\rfloor}, \Pi
\bigl(p_{\ov\be}(t)\bigr)\bigr)\dvtx  t\in[0,1]\bigr\}
\]
where $\Pi$ stands for the canonical projection from $\t_{\ov\be}$
onto $\bm_\infty$. The distortion of this
correspondence is
\begin{eqnarray*}
&&\sup_{0\leq s,t\leq1} \biggl| \biggl(\frac{3}{4}\biggr)^{3/8}n^{-1/4}{\mathrm{d}}^{\bq_n}_{\mathrm{gr}}\bigl(u^n_{\lfloor2ns\rfloor},
u^n_{\lfloor2nt\rfloor}\bigr) -D^*\bigl(p_{\ov\be}(s),p_{\ov\be}(t)
\bigr) \biggr|
\\
&&\quad= \sup_{0\leq s,t\leq1} \biggl| \biggl(\frac{3}{4}
\biggr)^{3/8}n^{-1/4} {\mathrm{d}}_n\bigl(\lfloor2ns
\rfloor,\lfloor2nt\rfloor\bigr) -D(s,t) \biggr|
\end{eqnarray*}
using Lemma~\ref{crucial-main} to write $D^*(p_{\ov\be}(s),p_{\ov
\be}(t))=D(p_{\ov\be}(s),p_{\ov\be}(t))=D(s,t)$.
From the (almost sure) convergence (\ref{cd-dist}), the quantity in
the last display tends to $0$
as $n\to\infty$ along the sequence $(n_k)_{k\geq1}$. From the
expression of the
Gromov--Hausdorff distance in terms of correspondences~\cite{BBI}, Theorem
7.3.25, we conclude that
\[
\bigl(V(\bq_n)\setminus\{\partial\}, \bigl(\tfrac{3}{4}
\bigr)^{3/8}n^{-1/4} {\mathrm{d}}^{\bq_n}_{\mathrm{gr}}
\bigr) \build{\la} {} {\mathrm{a.s.}} \bigl(\bm_\infty,D^*\bigr)
\]
as $n\to\infty$ along the sequence $(n_k)_{k\geq1}$. Clearly, the
latter convergence also
holds if we replace $V(\bq_n)\setminus\{\partial\}$ by $V(\bq_n)$.

The preceding arguments show that from any sequence of integers
converging to $\infty$ we can extract
a subsequence along which the convergence of Theorem~\ref{MainTh}
holds. This is enough
to prove Theorem~\ref{MainTh}.

\begin{pf*}{Proof of Lemma~\ref{crucial-main}}
Here we follow the ideas of the treatment of triangulations in
\cite{LGU}, Section~8. By a continuity argument, it is enough to prove that
if $X$ and $Y$ are independent and uniformly distributed over $[0,1]$,
and independent of the sequence $(\bq_n)_{n\geq1}$ (and therefore also
of $(\be,Z,D)$), we have
\[
D\bigl(p_{\ov\be}(X),p_{\ov\be}(Y)\bigr) = D^*\bigl(p_{\ov\be}(X),p_{\ov
\be
}(Y)
\bigr)\quad\mbox{a.s.}
\]
As we already know that $D(p_{\ov\be}(X),p_{\ov\be}(Y)) \leq
D^*(p_{\ov\be}(X),p_{\ov\be}(Y))$, it
will be sufficient to prove that these two random variables have the
same distribution. The distribution
of $D^*(p_{\ov\be}(X),p_{\ov\be}(Y))$ is identified in Corollary
7.3 of~\cite{LGU}:
\[
D^*\bigl(p_{\ov\be}(X),p_{\ov\be}(Y)\bigr)\build{=} {} {(\mathrm{d})}
Z_X - \min Z \build{=} {} {(\mathrm{d})} \ov Z_X.
\]

On the other hand, we can also derive the distribution of
$D(p_{\ov\be}(X),p_{\ov\be}(Y))=D(X,Y)$. For every $n\geq1$, we set
\[
i_n=\lfloor2nX\rfloor,\quad j_n=\lfloor2nY\rfloor,
\]
so that $i_n$ and $j_n$ are independent (and independent of $\bq_n$)
and uniformly distributed over $\{0,1,\ldots,2n-1\}$. Recall that
every integer
$i\in\{0,1,\ldots,2n-1\}$ corresponds to a corner of the tree $\tau
_n$ and
therefore via Schaeffer's bijection to an edge of $\bq_n$. We define
$\wt\bq_n$ by saying that $\wt\bq_n$ is the same planar map as $\bq_n$
but re-rooted at the edge associated with $i_n$, with each of the two possible
orientations chosen with probability $\frac{1}{2}$. Then
$\wt\bq_n$ is also uniformly distributed over $\mc Q_n^{\mathrm{nice}}$,
and we
let $\wt\tau_n$ be the associated tree in Schaeffer's bijection.
Write $\wt{\mathrm{d}}_n$
for the analog of ${\mathrm{d}}_n$ when $\bq_n$ is replaced by $\wt\bq_n$.

Let $k_n\in\{0,1,\ldots,2n-1\}$ be the index of the corner of the
tree $\wt\tau_n$
corresponding to the edge of $\bq_n$ that starts from the corner $j_n$
of $\tau_n$ in Schaeffer's bijection. Note that, conditionally on the pair
$(\bq_n,\wt\bq_n)$, the latter edge is uniformly distributed over
all edges
of $\bq_n$, and is thus also uniformly distributed over all edges of
$\wt\bq_n$
(recall that $\wt\bq_n$ is the same quadrangulation as $\bq_n$
with a different root). Hence, conditionally on the pair
$(\bq_n,\wt\bq_n)$, $k_n$ is uniformly distributed over $\{
0,1,\ldots,2n-1\}$, and in
particular the random variable $k_n$ is independent of $\wt\bq_n$.
We next observe that
%
\begin{equation}
\label{crucialbd}\bigl| {\mathrm{d}}_n(i_n,j_n) - \wt{\mathrm{d}}_n(0,k_n)\bigr| \leq2,
\end{equation}
because, with an obvious notation, the vertex $u^n_{i_n}$ is either
equal or adjacent to $\wt u^n_{0}$, and
similarly $u^n_{j_n}$ is either equal or adjacent to $\wt u^n_{k_n}$.

Now we have $\wt{\mathrm{d}}_n(0,k_n)\build{=}{}{(\mathrm{d})} {\mathrm{d}}_n(0,i_n)$ and by (\ref{cd-dist}),
\[
\bigl(\tfrac{3}{4}\bigr)^{3/8} n^{-1/4}{\mathrm{d}}_n(0,i_n)\build{\la} {n\to\infty}
{\mathrm{a.s.}} D(0,X)=\ov
Z_X
\]
where the convergence holds a.s. along the sequence $(n_k)_{k\geq1}$.
Similarly,
(\ref{cd-dist}) implies that, along the same sequence,
\[
\bigl(\tfrac{3}{4}\bigr)^{3/8} n^{-1/4}{\mathrm{d}}_n(i_n,j_n)\build{\la} {n\to\infty}
{\mathrm{a.s.}} D(X,Y).
\]
From the last two convergences and (\ref{crucialbd}), we obtain that
$D(X,Y)$ has the same distribution as~$\ov Z_X$. Since we already
observed that this
is also the distribution of $D^*(p_{\ov\be}(X),p_{\ov\be}(Y))$, the
proof of the lemma
is complete.
\end{pf*}

\section{The convergence of coding functions}\label{sec4}

In this section, we prove our main technical result Theorem~\ref{conv-coding}.
We start by deriving an intermediate convergence theorem.

\subsection{A preliminary convergence}\label{sec4.1}
\label{subsec-preli}

If $\tau$ is a plane tree, we let $\partial\tau$ stand for the set
of all leaves of $\tau$
different from the root vertex (which may or may not be a leaf). Then,
for every
integer $n\geq0$, we define $\mc W_n^{\circ}$ as the set of all labeled
trees $(\tau,U)$ such that $|\tau|=n$ and the property
$U(e_+)=U(e_-)-1$ holds for every
edge $e\in E(\tau)$ such that $e_+\in\partial\tau$. We also set
\[
\mc W^{\circ} = \bigcup_{n=0}^\infty
\mc W_n^{\circ}.
\]

Let $\beta:=\frac{1}{2}(\sqrt{3}-1)$ and let $\mu$ be the
probability measure on $\{0,1,2,\ldots\}$
defined by
\[
\mu(0):=\frac{1}{z_\beta} \frac{1}{3},\quad \mu(k):=\frac
{1}{z_\beta}
\beta^k\quad\mbox{for every }k\geq1,
\]
where $z_\beta$ is the appropriate normalizing constant:
\[
z_\beta= \frac{1}{3} + \sum_{k=1}^\infty
\beta^k = \frac{\sqrt{3}
+ 1}{3}.
\]
An easy calculation shows that $\mu$ is critical, meaning that
\[
\sum_{k=0}^\infty k\mu(k)= 1.
\]
In fact, the value of $\beta$ has been chosen so that this criticality
property holds. We can
also compute the variance of $\mu$,
\[
\sigma^2:= \sum_{k=0}^\infty(k-1)^2
\mu(k) = \frac{2}{\sqrt{3}}.
\]

Next, let $\t$ be a Galton--Watson tree with offspring distribution
$\mu$. Since
$\mu$ is critical, $\t$ is almost surely finite, and we can view $\t
$ as a
random variable with values in the space of all plane trees. We then
define random labels
$\wt\u(v)$, $v\in\t$ in the following way. We set $\wt\u
(\varnothing)=0$ and conditionally
on $\t$, we choose the other labels $\wt\u(v)$, $v\in\t\setminus
\{\varnothing\}$ in such a way
that the random variables $\wt\u(e_+)-\wt\u(e_-)$, $e\in E(\t)$,
are independent and uniformly
distributed over $\{-1,0,1\}$. In this way, we obtain a (random)
labeled tree
$(\t,\widetilde\u)$, and we may assume that $(\t,\widetilde\u)$
is also defined
on the probability space $(\Omega,\mathcal{A},\bb P)$.

There is of course no reason why the labeled tree $(\t,\widetilde\u
)$ should
belong to $\mc W^{\circ}$, and so we modify it in the following way.
We set $\u(v)=\wt\u(v)$
for every vertex $v\in\t\setminus\partial\t$. On the other hand,
for every edge
$e\in E(\t)$ such that $e_+\in\partial\t$, we set
$\u(e_+)=\wt\u(e_-) - 1$.
Then $(\t,\u)$ is a random element of $\mc W^{\circ}$.

The motivation for the preceding construction comes from the following lemma.

\begin{lemma}
\label{unifdi}
The conditional distribution of $(\t,\u)$ knowing that $|\t|=n$ is
the uniform
probability measure on $\mc W^{\circ}_n$.
\end{lemma}

\begin{pf} The case $n=0$ is trivial and we exclude it in the following
argument.
Let \mbox{$(\tau,U)\in\mc W^{\circ}_n$}. We have
\[
{\bb P}\bigl((\t,\u)=(\tau,U)\bigr)= {\bb P}(\t=\tau)\times\bigl(
\tfrac
{1}{3}\bigr)^{|\tau|-\#\partial\tau}
\]
since $|\tau|-\#\partial\tau$ is the number of edges $e\in E(\tau)$
such that
$e_+\notin\partial\tau$. On the other hand,
\[
{\bb P}(\t=\tau)= \prod_{v\in\tau} \mu
\bigl(k_v(\tau)\bigr)= \biggl(\frac
{1}{z_\beta}\biggr)^{|\tau|+1}
\times\biggl(\frac{1}{3}\biggr)^{\#\partial\tau}\times\prod
_{v\in\tau
\setminus\partial\tau} \beta^{k_v(\tau)}.
\]
Since $\sum_{v\in\tau\setminus\partial\tau}k_v(\tau)=|\tau
|=n$, we arrive at
\[
{\bb P}\bigl((\t,\u)=(\tau,U)\bigr)= \biggl(\biggl(\frac{1}{z_\beta}
\biggr)^{n+1} \times\biggl(\frac{1}{3}\biggr)^{\#\partial\tau}\times
\beta^n \biggr) \times\biggl(\frac{1}{3}\biggr)^{n-\#\partial\tau} =
\biggl(\frac{1}{z_\beta
}\biggr)^{n+1}\times\biggl(\frac{\beta}{3}
\biggr)^n.
\]
This quantity does not depend on the choice of $(\tau,U)\in\mc
W^{\circ}_n$, and
the statement of the lemma follows.
\end{pf}

We write $\mc C=(\mc C_t)_{t\geq0}$ and $\mc V=(\mc V_t)_{t\geq0}$
for the
contour function and the label function of the labeled tree $(\t,\u
)$. We define
rescaled versions of $\mc C$ and $\mc V$ by setting for every $n\geq1$
and $t\in[0,1]$,
%
\begin{equation}
\label{scalingcontour} \mc C^n_t:= \frac{\sigma}{2}
n^{-1/2} \mc C_{2nt},\quad \mc V^n_t=
\biggl(\frac{2}{3}\biggr)^{-1/2} \biggl(\frac{\sigma}{2}
\biggr)^{1/2}n^{-1/4} \mc V_{2nt}
\end{equation}
where we recall that $\sigma^2= 2/\sqrt{3}$ is the variance of $\mu$
(in the previous display, $\frac{2}{3}$
corresponds to the variance of the uniform distribution on $\{-1,0,1\}
$). Note that
\[
\frac{\sigma}{2}=12^{-1/4},\quad\biggl(\frac{2}{3}
\biggr)^{-1/2} \biggl(\frac{\sigma}{2}\biggr)^{1/2}= \biggl(
\frac{3}{4}\biggr)^{3/8}.
\]

We write ${\bb P}_n={\bb P}( \cdot\mid|\t|=n)$ for the
conditional probability knowing that $|\t|=n$, and ${\bb E}_n$ for the
expectation under ${\bb P}_n$.

\begin{proposition}
\label{prelimiconv}
The law of $(\mc C^n_t,\mc V^n_t)_{0\leq t\leq1}$ under ${\bb P}_n$
converges as
$n\to\infty$ to the distribution of $(\be_t,Z_t)_{0\leq t\leq1}$.
\end{proposition}

\begin{pf}
Let $\wt{\mc V}$ stand for the label function of the labeled tree
$(\t,\wt\u)$. By construction, we have $|\u(v)-\wt\u(v)|\leq2$ for
every $v\in\t$, and it follows that for every $t\geq0$,
\[
|\wt{\mc V}_t -\mc V_t| \leq2.
\]
Let $\wt{\mc V}^n$ be defined from $\wt{\mc V}$ by the same scaling
operation we used to
define $\mc V^n$ from $\mc V$. From the preceding bound, we have also,
for every
$t\geq0$,
%
\begin{equation}
\label{boundpreli} \bigl|\wt{\mc V}^n_t -\mc
V^n_t\bigr| \leq\sqrt{3} \sigma^{1/2}
n^{-1/4}.
\end{equation}
By known results about the convergence of
discrete snakes~\cite{JS} (see Theorem 2.1 in~\cite{LG0}), we know
that the law
of $(\mc C^n_t,\wt{\mc V}^n_t)_{0\leq t\leq1}$ under
${\bb P}_n$ converges as
$n\to\infty$ to the distribution of $(\be_t,Z_t)_{0\leq t\leq1}$.
The statement of the proposition immediately follows from this
convergence and
the bound (\ref{boundpreli}).
\end{pf}

We will be interested in conditional versions of the convergence of
Proposition~\ref{prelimiconv}. Let us start by discussing a simple
case. For every
real $x\geq0$, we write $\ov{\bb P}{}^x_n$ for the conditional
probability measure
\[
\ov{\bb P}{}^x_n= {\bb P}_n\bigl( \cdot\mid
\u(v)\geq-x\mbox{ for every }v\in\t\bigr).
\]
We write\vspace*{1pt} $\ov{\bb P}_n=\ov{\bb P}{}^0_n$ to simplify notation. We
denote the expectation
under $\ov{\bb P}{}^x_n$, respectively, under $\ov{\bb P}_n$, by $\ov{\bb
E}{}^x_n$, respectively,
$\ov{\bb E}_n$.

Let $(r_n)_{n\geq1}$ be a sequence of positive real numbers converging
to $r>0$,
and let $F$ be a bounded continuous function on $C([0,1],\bb R^2)$. It
follows from
the preceding proposition (together with the fact that the law of
$\min_{0\leq t\leq1} Z_t$ has no atoms) that
\[
\lim_{n\to\infty} {\bb E}_n\bigl[F\bigl(\mc
C^n,\mc V^n\bigr) {\mathbf1}_{\{\min
_{0\leq t\leq1}\mc V^n_t \geq-r_n\}}\bigr] = {\bb E}
\bigl[F(\be,Z) {\mathbf1}_{\{\min_{0\leq t\leq1} Z_t\geq-r\}}\bigr].
\]
Let $\kappa_n:= \frac{2}{\sqrt{3}}\sigma^{-1/2} n^{1/4}$
be the inverse of the scaling factor in the definition of $\mc V^n$.
The preceding convergence
implies that
\[
\lim_{n\to\infty} \ov{\bb E}{}^{\kappa_nr_n}_n\bigl[F
\bigl(\mc C^n,\mc V^n\bigr)\bigr] = {\bb E}\bigl[F\bigl(
\be^{(r)},Z^{(r)}\bigr)\bigr].
\]
Since this holds for any sequence $(r_n)$ converging to $r>0$, we get
that, for any compact
subinterval $I$ of $(0,\infty)$, we have also
%
\begin{equation}
\label{consnakecondi} \lim_{n\to\infty} \sup_{r\in I}
\bigl|\ov{\bb E}{}^{\kappa_nr}_n\bigl[F\bigl(\mc C^n,\mc
V^n\bigr)\bigr] - {\bb E}\bigl[F\bigl(\be^{(r)},Z^{(r)}
\bigr)\bigr] \bigr| = 0.
\end{equation}

A labeled tree codes a nice quadrangulation with $n$ faces if and only
if it is a
tree of $\mc W_n^{\circ}$ with nonnegative labels,
and, in the case when the root is a leaf, if the label of the only
child of the root is $1$ and if there is
another vertex with label $0$. Recalling
Lemma~\ref{unifdi}, we see that scaling limits for the contour and
label functions of a labeled tree uniformly distributed over $\mc
W_n^{\circ}$
are given by the preceding proposition. As in the previous discussion,
Theorem~\ref{conv-coding}
can thus be seen as a conditional version of Proposition \ref
{prelimiconv}. Closely related
conditionings are discussed in~\cite{LG0}, but we shall not be able
to apply directly the results of~\cite{LG0} (though we use certain
ideas of
the latter paper).

Let $H$ be the set of all labeled trees $(\tau,U)$
such that:
\begin{enumerate}
\item[$\bullet$] either $k_\varnothing(\tau)\geq2$;
\item[$\bullet$] or $k_\varnothing(\tau)=1$, $U(1)=1$, and there exists
$v\in\tau\setminus\{\varnothing\}$ such that $U(v)=0$.
\end{enumerate}
By Proposition~\ref{nicetree}, the set of all labeled trees associated
with nice quadrangulations with $n$ faces (in
Schaeffer's bijection) is
%
\begin{equation}
\label{treenice} \mc W^{\mathrm{nice}}_n = \mc W_n^\circ
\cap\mc W_n^+ \cap H.
\end{equation}
We write $\mc H$ for the event $\mc H:=\{(\t,\u)\in H\}$.

\subsection{A spatial Markov property}\label{sec4.2}
We consider again the random labeled tree $(\t,\u)$ introduced in the
previous subsection.
A major difficulty in the proof of Theorem~\ref{conv-coding} comes
from the fact
that conditioning the tree on having nonnegative labels is not easy to handle.
To remedy this problem, we will introduce a (large) subtree of $\t$, which
in a sense will approximate $\t$, but whose distribution will involve
a less
degenerate conditioning (see Proposition~\ref{keytech} below).

Recall the notation $\bb V,\bb V^{(w)},\bb V_{(w)}$
introduced in Section~\ref{sec2}. Let $w\in\bb V$, and first argue
on the event $\{w\in\t\}$. We let
\[
\t^{(w)}:=\t\cap\bb V^{(w)}
\]
be the set of all vertices of $\t$ that are not strict descendants of $w$.
Clearly, $\t^{(w)}$ is a tree and we equip it with labels by setting
$\u^{(w)}(v)=\u(v)$ for every $v\in\t^{(w)}$.
We similarly\vspace*{1pt} define
$\wt\u^{(w)}(v)=\wt\u(v) $ for every $v\in\t^{(w)}$. If $w\notin
\t$, we just
put $\t^{(w)}=\{\varnothing\}$ and $\wt\u^{(w)}(\varnothing)=\u
^{(w)}(\varnothing)=0$.

Next, on the event $\{w\in\t\}$, we define
\[
\t_{(w)}=\{v\in\bb V\dvtx  wv\in\t\}.
\]
Then $\t_{(w)}$ is a tree (we may view it as the subtree of
descendants of $w$). We assign
labels to the vertices of $\t_{(w)}$ by setting, for $v\in\t_{(w)}$,
\[
\u_{(w)}(v)=\u(wv)-\u(w).
\]
On the event $\{w\notin\t\}$ we set
$\t_{(w)}=\{\varnothing\}$ and $\u_{(w)}(\varnothing)=0$.

For every $w\in\bb V$,
let $\Sigma^{(w)}$ be the $\sigma$-field generated by $(\t^{(w)},\wt
\u^{(w)})$.

\begin{lemma}
\label{lem-spa}
For every nonnegative function $G$ on the space of all labeled trees, for
every $w\in\bb V$,
\[
{\bb E}\bigl[{\mathbf1}_{\{w\in\t\}}G(\t_{(w)},\u_{(w)})\mid
\Sigma^{(w)}\bigr] = {\mathbf1}_{\{w\in\t\}} {\bb E}\bigl[G(\t,\u)\bigr].
\]
\end{lemma}

\begin{remark} It is essential that we define $\Sigma^{(w)}$ as the
$\sigma$-field generated by
the pair $(\t^{(w)},\wt\u^{(w)})$, and not by the pair $(\t
^{(w)},\u^{(w)})$: The knowledge of $(\t^{(w)},\u^{(w)})$ provides
information about the fact that $w$ is or is not a leaf of $\t$, and
the statement of the
lemma would not hold with this alternative definition.
\end{remark}

The proof of Lemma~\ref{lem-spa} is a simple application of properties of
Galton--Watson trees and the way labels are generated. We omit the details.

Let us introduce some notation.
We fix an integer $x\geq1$ and define a subset $\ell_x$ of $\t$ by setting
\[
\ell_{x} = \bigl\{ w\in\mc T\dvtx  \mc U(w) \ge x \mbox{ and } \mc
U(v)<x \mbox{ for every } v \in[\varnothing, w]\setminus\{w\}\bigr\},
\]
where
$
[\varnothing,w]=\{v\in\bb V\dvtx  w\in\bb V_{(v)}\} $
stands for the set of all ancestors of $w$. Define $\wt\ell^x$
similarly by replacing $\u$ by $\wt\u$.

Next fix $r\in(1/2,1)$ and for every $n\geq1$, consider the event
\[
F^{r,x}_n:=\bigl\{ |\mc T|=n \mbox{ and there exists } w\in\ell
_{x} \mbox{ such that } |\mc T_{(w)}| \ge rn\bigr\}.
\]
If $F^{r,x}_n$ holds, the vertex $w\in\ell_{x}$ such that $|\mc
T_{(w)}| \ge rn$ is clearly unique, and
we denote it by $w_n$. We also set $m_n=|\t_{(w_n)}|=n-|\t^{(w_n)}|$
on the same event. If
$F^{r,x}_n$ does not hold, we set $w_n=\varnothing$ and $m_n=n$ for
definiteness.

The following technical result plays a major role in our
proof of Theorem~\ref{conv-coding}. Roughly speaking,
this result identifies the distribution, under the probability measure
$\ov{\bb P}_n$ restricted to the event~$F^{r,x}_n$,
of the ``large'' subtree of $\mc T$ rooted at the vertex $w_n$.

\begin{proposition}
\label{keytech}
Let $G_1,G_2$ be nonnegative functions on the space $\mc W$. Then,
\begin{eqnarray*}
&&\ov{\bb E}_n \bigl[ {\mathbf1}_{F^{r,x}_n} G_1\bigl(
\t^{(w_n)},\u^{(w_n)}\bigr) G_2(\t_{(w_n)},
\u_{(w_n)}) \bigr]
\\
&&\quad= \ov{\bb E}_n \bigl[ {\mathbf1}_{F^{r,x}_n} G_1
\bigl(\t^{(w_n)},\u^{(w_n)}\bigr) \ov{\bb E}{}^x_{m_n}
\bigl[G_2(\t,\u)\bigr] \bigr].
\end{eqnarray*}
\end{proposition}

\begin{pf}
We fix $w\in\bb V \setminus\{\varnothing\}$ and $m\in[rn,n]\cap
\bb Z$. On the event $\{w\in\t\}$,
we also
set $\ov\u{}^{(w)}(v)=\u^{(w)}(v)$ if $v\in\t^{(w)}\setminus\{w\}$
and $\ov\u{}^{(w)}(w)=\wt\u^{(w)}(w)$.
Then the quantity
\[
{\mathbf1}_{F^{r,x}_n\cap\{\u(v)\geq0,\forall v\in\t\}\cap\{
(w_n,m_n)=(w,m)\}} G_1\bigl(\t^{(w)},\u^{(w)}
\bigr) G_2(\t_{(w)},\u_{(w)})
\]
is equal to the product of
\[
R:= {\mathbf1}_{\{w\in\wt\ell_x\}\cap
\{|\t^{(w)}|=n-m\}\cap\{\u(v)\geq0,\forall v\in\t^{(w)}\setminus
\{w\}\}} G_1\bigl(\t^{(w)},\ov
\u{}^{(w)}\bigr)
\]
with
\[
S:= {\mathbf1}_{\{|\t_{(w)}|=m\}\cap\{\u(v)\geq0,\forall v\in\t\cap
\bb V_{(w)}\}} G_2(\t_{(w)},\u_{(w)}).
\]
The point is that if $|\t_{(w)}|=m>0$, $w$ is a vertex of $\t$ that
is not a leaf, so that $\wt\u^{(w)}(w)=\u^{(w)}(w)$,
implying that $G_1(\t^{(w)},\ov\u{}^{(w)})= G_1(\t^{(w)},\u^{(w)})$
and that the property $w\in\ell_x$
holds if and only if $w\in\wt\ell_x$.

It is easy to verify that $R$ is $\Sigma^{(w)}$-measurable. Notice that
${\mathbf1}_{\{w\in\t\}} G_1(\t^{(w)},\ov\u{}^{(w)})$ is $\Sigma
^{(w)}$-measurable, which would not
be the case for ${\mathbf1}_{\{w\in\t\}} G_1(\t^{(w)},\u^{(w)})$.

Next, notice that on the event $\{w\in\wt\ell_x\}\cap\{|\t
_{(w)}|=m\}
= \{w\in\ell_x\}\cap\{|\t_{(w)}|=m\}$, we have necessarily $\u(w)=x$
and the property
\[
\u(v)\geq0\quad\forall v\in\t\cap\bb V_{(w)}
\]
holds if and only if
\[
\u_{(w)}(v)\geq-x\quad \forall v\in\t_{(w)}.
\]
This shows that, on the event $\{w\in\wt\ell_x\}$, $S$ coincides
with the variable
\[
{\mathbf1}_{\{|\t_{(w)}|=m\}\cap\{\u_{(w)}(v)\geq-x, \forall v\in\t
_{(w)}\}} G_2(\t_{(w)},\u_{(w)})
\]
which is a function of the pair $(\t_{(w)},\u_{(w)})$.

Recalling that $R$ is $\Sigma^{(w)}$-measurable and using Lemma \ref
{lem-spa}, we get
\begin{eqnarray*}
{\bb E}[RS]&=& {\bb E} \bigl[R {\mathbf1}_{\{|\t_{(w)}|=m\}\cap\{\u
_{(w)}(v)\geq-x, \forall v\in\t_{(w)}\}} G_2(
\t_{(w)},\u_{(w)}) \bigr]
\\
&=&{\bb E}[R]\times{\bb E} \bigl[{\mathbf1}_{\{|\t|=m\}\cap\{\u(v)\geq
-x,\forall v\in\t\}} G_2(\t,\u)
\bigr]
\\
&=&{\bb E}[R]\times{\bb P} \bigl[\bigl\{|\t|=m\bigr\}\cap\bigl\{\u(v)\geq
-x,\forall v\in
\t\bigr\} \bigr]\times\ov{\bb E}{}^x_m \bigl[G_2(
\t,\u) \bigr]
\end{eqnarray*}
by the definition of the conditional measure $\ov{\bb P}{}^x_m$.

From the case $G_2=1$ in the equality between the two ends of the last
display, we have also
\begin{eqnarray*}
&&{\bb E}[R]\times{\bb P} \bigl[\bigl\{|\t|=m\bigr\}\cap\bigl\{\u(v)\geq
-x,\forall v\in\t
\bigr\} \bigr]
\\
&&\quad= {\bb E} \bigl[{\mathbf1}_{F^{r,x}_n\cap\{\u(v)\geq0,\forall
v\in\t\}\cap\{(w_n,m_n)=(w,m)\}} G_1\bigl(
\t^{(w)},\u^{(w)}\bigr) \bigr].
\end{eqnarray*}
By substituting this into the preceding display, we arrive at
\begin{eqnarray*}
&&{\bb E} \bigl[ {\mathbf1}_{F^{r,x}_n\cap\{\u(v)\geq0,\forall v\in\t\}
\cap\{(w_n,m_n)=(w,m)\}} G_1\bigl(\t^{(w)},
\u^{(w)}\bigr) G_2(\t_{(w)},\u_{(w)})
\bigr]
\\
&&\quad ={\bb E} \bigl[{\mathbf1}_{F^{r,x}_n\cap\{\u(v)\geq0,\forall
v\in\t\}\cap\{(w_n,m_n)=(w,m)\}} G_1\bigl(\t^{(w)},
\u^{(w)}\bigr) \bigr]\times\ov{\bb E}{}^x_m
\bigl[G_2(\t,\u) \bigr]
\\
&&\quad ={\bb E} \bigl[{\mathbf1}_{F^{r,x}_n\cap\{\u(v)\geq0,\forall
v\in\t\}\cap\{(w_n,m_n)=(w,m)\}} G_1\bigl(\t^{(w)},
\u^{(w)}\bigr) \ov{\bb E}{}^x_{m_n}
\bigl[G_2(\t,\u) \bigr] \bigr].
\end{eqnarray*}
Now we just have to sum over all possible choices of $w$ and $m$ and
divide by the quantity
${\bb P}(\{|\t|=n\}\cap\{\u(v)\geq0,\forall v\in\t\})$ to get the
statement of the proposition.
\end{pf}

\subsection{Technical estimates}\label{sec4.3}

Recall from Section~\ref{subsec-preli} the definition of the set
$\mc H$ and
of the rescaled processes $\mc C^n$ and~$\mc V^n$. To simplify
notation, we write
$\ov{\bb P}_{n,\mc H}$ for the conditional probability $\ov{\bb
P}_n(\cdot\mid\mc H)$. This makes
sense as soon as $\ov{\bb P}_n(\mc H)>0$, which holds for every $n\geq2$.

To simplify notation, we write $k_v$ instead of $k_v(\t)$ in the following.

\begin{proposition}
\label{techest}
There exists a constant $a_0>0$ such that $\ov{\bb P}_n(\mc H)\geq
a_0$ for every $n\geq2$.
Moreover, for any $b>0$ and $\vep>0$, we can find $\delta,\alpha\in
(0,\frac{1}{4})$
such that, for every sufficiently large $n$,
\[
\ov{\bb P}_{n,\mc H} \Bigl( \inf_{t\in[\delta,1-\delta]} \mc
V^n_t > \alpha, \sup_{t\in
[0,2\delta]\cap[1-2\delta,1]} \bigl(\mc
C^n_t + \mc V^n_t\bigr)\leq\vep
\Bigr) \geq1-b.
\]
\end{proposition}

\begin{pf}
We start by proving the first assertion. It is enough to find a
constant $a_0$ such that, for every sufficiently large
$n$,
%
\begin{equation}
\label{claimt} {\bb P}_n\bigl(k_\varnothing\geq2; \u(v)\geq0,
\forall v\in\t\bigr) \geq a_0 {\bb P}_n\bigl(\u(v)\geq0,
\forall v\in\t\bigr).
\end{equation}
Now observe that, by construction,
%
\begin{equation}
\label{trivialbounds} \min_{v\in\t} \wt{\mc U}(v) \geq\min
_{v\in\t} \mc U(v) \geq\min_{v\in\t} \wt{\mc
U}(v) -1.
\end{equation}
In particular, it is immediate that
%
\begin{equation}
\label{claimt1} {\bb P}_n\bigl(\u(v)\geq0,\forall v\in\t\bigr)\leq{
\bb P}_n\bigl(\wt\u(v)\geq0,\forall v\in\t\bigr).
\end{equation}
On the other hand, we get a lower bound on ${\bb P}_n(k_\varnothing
\geq2;\u(v)\geq0,\forall v\in\t)$ by considering the event
where $\varnothing$ has (exactly) two children, who both have label
$1$, and the second child
of $\varnothing$ has one child, and this child is a leaf.
We get, for $n\geq4$,
%
\begin{eqnarray}
\label{claimt2} &&{\bb P}_n\bigl(k_\varnothing\geq2; \u(v)\geq0,
\forall v\in\t\bigr)
\nonumber
\\
&&\quad = \frac{{\bb P}(k_\varnothing\geq2; \u(v)\geq0,\forall
v\in\t;|\t|=n)}{{\bb P}(|\t|=n)}
\nonumber
\\
&&\quad \geq\frac{{\bb P}(k_\varnothing= 2,k_2=1,k_{21}=0,\u(1)=\u
(2)=1; \u(v)\geq0,\forall v\in\t;|\t|=n)}{{\bb P}(|\t|=n)}
\nonumber\\[-8pt]\\[-8pt]
&&\quad = \frac{1}{{\bb P}(|\t|=n)} \frac{\beta^2}{z_\beta
}\times\frac{\beta}{z_\beta} \times
\frac{1}{3z_\beta} \times\biggl(\frac{1}{3}\biggr)^2\times{\bb
P}\bigl(\u(v)\geq-1,\forall v\in\t;|\t|=n-3\bigr)\quad
\nonumber
\\
&&\quad \geq\frac{c}{{\bb P}(|\t|=n)} {\bb P}\bigl(\wt\u(v)\geq0,\forall
v\in\t;|\t|=n-3
\bigr)
\nonumber
\\
&&\quad = c \frac{{\bb P}(|\t|=n-3)}{{\bb P}(|\t|=n)} {\bb P}_{n-3}\bigl(\wt
\u(v)\geq0,\forall v
\in\t\bigr)\nonumber
\end{eqnarray}
where $c= \beta^3/(27z_\beta^3)$.
Proposition 4.2 in~\cite{LG0} gives the existence of two positive constants
$c_1$ and $c_2$
such that, for every sufficiently large $n$,
%
\begin{equation}
\label{boundposi} \frac{c_1}{n} \leq{\bb P}_{n}\bigl(\wt\u(v)
\geq0,\forall v\in\t\bigr)\leq\frac{c_2}{n}.
\end{equation}
Moreover, standard asymptotics for the total progeny of a critical Galton--Watson
tree show that
\[
\lim_{n\to\infty} \frac{{\bb P}(|\t|=n-3)}{{\bb P}(|\t|=n)} =1.
\]
Our claim (\ref{claimt}) follows from the preceding observations
together with the
bounds (\ref{claimt1}) and (\ref{claimt2}).

Let us turn to the proof of the second assertion. We start by observing
that the first
part of the proof, and in particular (\ref{claimt2}) and (\ref
{boundposi}) show that the bounds
%
\begin{equation}
\label{boundposi2} \frac{c'_1}{n} \leq{\bb P}_{n} \Bigl(\min
_{v\in\t} \u(v)\geq0 \Bigr)\leq{\bb P}_{n} \Bigl(\min
_{v\in\t} \wt\u(v)\geq0 \Bigr)\leq\frac{c_2}{n}
\end{equation}
hold for every sufficiently large $n$, with a positive constant $c'_1$.
Then, thanks to the lower bound $\ov{\bb P}_n(\mc H)\geq a_0$, it is
enough to verify that
given $b>0$ and $\vep>0$, we can find $\delta,\alpha\in(0,\vep
\wedge\frac{1}{4})$
so that
\[
\ov{\bb P}_{n} \Bigl( \Bigl\{ \inf_{t\in[\delta,1-\delta]} \mc
V^n_t \leq\alpha\Bigr\}\cup\Bigl\{ \sup
_{t\in[0,2\delta]\cap[1-2\delta,1]} \bigl(\mc C^n_t + \mc
V^n_t\bigr)> \vep\Bigr\} \Bigr) \leq b.
\]
Recall the notation $\wt{\mc V}^n$ introduced in the proof of
Proposition~\ref{prelimiconv}
and the bound (\ref{boundpreli}). Clearly, it is enough to verify that
the bound of the preceding
display holds when $\mc V^n$ is replaced by $\wt{\mc V}^n$. To
simplify notation, set
\[
A_n^{\delta,\alpha}= \Bigl\{ \inf_{t\in[\delta,1-\delta]}\wt{ \mc
V}^n_t \leq\alpha\Bigr\} \cup\Bigl\{ \sup
_{t\in[0,2\delta]\cap[1-2\delta,1]} \bigl(\mc C^n_t + \wt{\mc
V}^n_t\bigr)> \vep\Bigr\}.
\]
We have then
\begin{eqnarray*}
\ov{\bb P}_n\bigl(A_n^{\delta,\alpha}\bigr) &=&
\frac{{\bb P}_n(A_n^{\delta,\alpha}\cap\{\min_{v\in\t} \u
(v)\geq0\})}{{\bb P}_n(\min_{v\in\t} \u(v)\geq0)}
\\
&\leq& \frac{{\bb P}_n(\min_{v\in\t} \wt\u(v)\geq0)}{{\bb
P}_n(\min_{v\in\t} \u(v)\geq0)}\times{\bb P}_n \Bigl(A_n^{\delta,\alpha}
\bigm|
\min_{v\in\t} \wt\u(v)\geq0 \Bigr),
\end{eqnarray*}
using (\ref{trivialbounds}) in the last bound.
On the one hand, the bounds from (\ref{boundposi2}) imply that the ratio
\[
\frac{{\bb P}_n(\min_{v\in\t} \wt\u(v)\geq0)}{{\bb P}_n(\min_{v\in\t} \u
(v)\geq0)}
\]
is bounded above by a constant. On the other hand Proposition 6.1 in
\cite{LG0}
shows that
the quantity
\[
{\bb P}_n \Bigl(A_n^{\delta,\alpha} \bigm| \min
_{v\in\t} \wt\u(v)\geq0 \Bigr)
\]
can be made arbitrarily small (for all sufficiently large $n$) by choosing
$\delta$ and $\alpha$ sufficiently small. This completes the proof of
the proposition.
\end{pf}

We write
\[
\Gamma^{\alpha,\delta,\vep}_n:=\bigl\{|\t|=n\bigr\}\cap\Bigl\{ \inf
_{t\in[\delta,1-\delta]} \mc V^n_t > \alpha, \sup
_{t\in
[0,2\delta]\cap[1-2\delta,1]} \bigl(\mc C^n_t + \mc
V^n_t\bigr)\leq\vep\Bigr\}
\]
for the event considered in Proposition~\ref{techest}.

Now recall the notation $F_{n}^{r,x}$ introduced before Proposition
\ref{keytech}.
Also recall the definition of the constants $\kappa_n$ a little before
(\ref{consnakecondi}). For
$\alpha>0$ and $\delta\in(0,\frac{1}{4})$, we set
\[
E_n^{\alpha,\delta} = F_n^{1-2\delta,\lfloor\alpha\kappa_n\rfloor}
\]
to simplify notation. We implicitly consider only values of $n$ such that
$ \alpha\kappa_n \geq1$. On the event $E_n^{\alpha,\delta}$, there
is a
unique vertex $w\in\ell_{\lfloor\alpha\kappa_n\rfloor}$ such that
$|\t_{(w)}|\geq(1-2\delta)n$
and we denote this vertex by $w_n^{\alpha,\delta}$ (as previously, if
$E_n^{\alpha,\delta}$
does not hold, we take $w_n^{\alpha,\delta}=\varnothing$). We also
set $m_n^{\alpha,\delta}= |\t_{(w_n^{\alpha,\delta})}|$.

\begin{lemma}
\label{lemtech0}
For every $\alpha>0$, $\delta\in(0,\frac{1}{4})$ and $\vep>0$,
we have $\Gamma^{\alpha,\delta,\vep}_n\subset E_n^{\alpha,\delta} $.
\end{lemma}

\begin{pf} Suppose that $\Gamma^{\alpha,\delta,\vep}_n$ holds. Then all
vertices of $\t$ visited by the contour exploration at integer times
between $2\delta n$
and $2(1-\delta)n$ must have a label strictly greater than $\alpha
\kappa_n$. By the properties
of the contour exploration, this implies that all these vertices share
a common ancestor $v_n$
belonging to $\ell_{\lfloor\alpha\kappa_n\rfloor}$, which moreover
is such that
$|\t_{(v_n)}|\geq(1-2\delta)n$. It follows that $E_n^{\alpha,\delta
}$ holds.
\end{pf}

We set
\[
\wt E_n^{\alpha,\delta}= E_n^{\alpha,\delta} \cap\bigl(\{
k_\varnothing\geq2\} \cup\bigl\{k_\varnothing=1,\u(1)=1,\exists v\in
\t^{(w_n^{\alpha,\delta
})}\setminus\{\varnothing\}\dvtx  \u(v)=0\bigr\} \bigr).
\]

\begin{proposition}
\label{techest2}
For any $b>0$, we can find $\delta,\alpha\in(0,\frac{1}{4})$
such that, for every sufficiently large~$n$,
\[
\ov{\bb P}_{n,\mc H} \bigl(\wt E^{\alpha,\delta}_n\bigr)
\geq1-b.
\]
\end{proposition}

\begin{pf}
By Proposition~\ref{techest} and Lemma~\ref{lemtech0}, it is enough to
verify that $\ov{\bb P}_{n,\mc H} (E^{\alpha,\delta}_n\setminus\wt
E^{\alpha,\delta}_n)$ tends to $0$ as $n\to\infty$, for any choice of
$\delta,\alpha\in(0,\frac{1}{4})$. By the first assertion of
Proposition~\ref{techest}, it suffices to verify that
%
\begin{equation}
\label{techest22} \lim_{n\to\infty} \ov{\bb P}_{n}
\bigl(\mc H\cap\bigl(E^{\alpha,\delta
}_n\setminus\wt
E^{\alpha,\delta}_n\bigr)\bigr)=0.
\end{equation}

Now observe that, on the event $\mc H\cap(E^{\alpha,\delta
}_n\setminus\wt E^{\alpha,\delta}_n)$, we have necessarily
$k_\varnothing=1$ and moreover there exists $v\in\t\setminus\t
^{(w_n^{\alpha,\delta})}$ such that $\u(v)=0$.
Consequently,
%
\begin{eqnarray}
\label{techest23} \ov{\bb P}_{n} \bigl(\mc H\cap
\bigl(E^{\alpha,\delta}_n\setminus\wt E^{\alpha,\delta}_n \bigr)\bigr)&
\leq&\ov{\bb P}_{n} \bigl(E^{\alpha,\delta}_n \cap\bigl\{\exists
v\in\t_{(w_n^{\alpha,\delta})}\dvtx \u(v)=0\bigr\}\bigr)
\nonumber\\[-8pt]\\[-8pt]
& = & \ov{\bb E}_n \bigl[ {\mathbf1}_{E^{\alpha,\delta}_n} \ov{\bb
P}{}^{\lfloor\alpha\kappa_n\rfloor}_{m_n^{\alpha,\delta
}}\bigl(\exists v\in\t\dvtx  \u(v)=-\lfloor\alpha
\kappa_n\rfloor\bigr) \bigr]\nonumber
\end{eqnarray}
using Proposition~\ref{keytech} in the last equality.

By construction, we have $n\geq m_n^{\alpha,\delta} \geq(1-2\delta
)n \geq n/2$ on the event $E^{\alpha,\delta}_n$.
An easy application of Proposition~\ref{prelimiconv} shows that
\[
\min_{\lceil n/2\rceil\leq m\leq n} {\bb P}_m \Bigl(\min
_{v\in\t} \u(v) \geq- \lfloor\alpha\kappa_n\rfloor
\Bigr) \geq c_{(\alpha)}
\]
with a constant $c_{(\alpha)}>0$ depending only on $\alpha$. Again
using Proposition~\ref{prelimiconv}
together with the fact that the law of $\inf_{t\in[0,1]} Z_t$ has no
atoms, we get that
\[
\sup_{\lceil n/2\rceil\leq m\leq n} {\bb P}_m \Bigl(\min
_{v\in\t} \u(v)=-\lfloor\alpha\kappa_n\rfloor\Bigr)
\build{\longrightarrow} {n\to\infty} {} 0.
\]
By combining the two preceding observations, we obtain that
\[
\sup_{\lceil n/2\rceil\leq m\leq n} \ov{\bb P}{}^{ \lfloor\alpha
\kappa_n\rfloor}_m \Bigl(\min
_{v\in\t} \u(v)= -\lfloor\alpha\kappa_n\rfloor\Bigr)
\build{\longrightarrow} {n\to\infty} {} 0.
\]
Our claim (\ref{techest22}) now follows from (\ref{techest23}).
\end{pf}

\subsection{Proof of the convergence of coding functions}\label{sec4.4}

We now turn to the proof of Theorem~\ref{conv-coding}.
Let us briefly discuss the main idea of the proof. We observe that, if
$\alpha$ and $\delta$ are
small enough, the
tree associated with a nice quadrangulation with $n$ faces is well
approximated by the subtree rooted at the vertex $w^{\alpha,\delta}_n$
introduced before Lemma~\ref{lemtech0}, whose label is small but
non-vanishing even
after rescaling. Together with Proposition~\ref{keytech}, the
convergence result (\ref{consnakecondi}) can then be used to relate the
law of this subtree and its labels to a conditioned pair $(\be
^{(r)},Z^{(r)})$. However, when $r$ is small
we know that the distribution of $(\be^{(r)},Z^{(r)})$ is close to
that of $(\be^{(0)},Z^{(0)})$.

We equip the space $C([0,1],\bb R^2)$ with the norm $\Vert(g,h)\Vert=
\Vert g \Vert_{\infty}\lor\Vert h\Vert_{\infty}$, where $\Vert
g\Vert_{\infty}$ stands for the supremum norm of $g$. For every $g\in
C([0,1],\bb R)$, and every $s>0$, we set:
\[
\omega_{g}(s) = \sup_{t_1,t_2\in[0,1], |t_1-t_2|\le s} \bigl|g({t_1})-g({t_2})\bigr|.
\]
We fix a Lipschitz function $F$ on $C([0,1],\bb R^2)$, with Lipschitz
constant less than $1$ and such that $0\le F\le1$.
By Lemma~\ref{unifdi} and (\ref{treenice}), the uniform distribution
on the space
$\mc W_n^{\mathrm{nice}}$ of all labeled trees asssociated with nice
quadrangulations with $n$ faces coincides with the law of $(\t,\u)$
under $\ov{\bb P}_{n,\mc H}$.
Therefore, to prove Theorem~\ref{conv-coding}, it is enough to show that
\[
\lim_{n\to\infty} \ov{\bb E}_{n,\mc H}\bigl[ F\bigl(\mc
C^n, \mc V^n\bigr)\bigr] = {\bb E} \bigl[F\bigl(
\be^{(0)},Z^{(0)}\bigr)\bigr].
\]
In the remaining part of this section we establish this convergence. To
this end, we fix $b>0$.

For every $\vep\in(0,\frac{1}{4})$ and $g\in C([0,1],\bb R)$, we set
\[
G_\vep(g)=\bigl( \omega_g(3\vep) +\bigl(4+2\| g
\|_\infty\bigr)\vep\bigr)\wedge1.
\]
For $r>0$, recall our notation $(\be^{(r)},Z^{(r)})$ for a process
whose distribution is the
conditional distribution of $(\be,Z)$ knowing that $\min_{0\leq t\leq
1}Z_t>-r$ (see the discussion
in Section~\ref{scalicofu}).
Since the distribution of $(\be^{(r)},Z^{(r)})$
depends continuously on $r\in[0,1]$, a simple argument shows that we
can choose
$\vep>0$ sufficiently small so that
%
\begin{equation}
\label{ini1} \sup_{r\in[0,1]} {\bb E}\bigl[G_\vep\bigl(
\be^{(r)}\bigr)+G_\vep\bigl(Z^{(r)}\bigr)\bigr] <
\frac{a_0b}{2},
\end{equation}
where we recall that the constant $a_0$ was introduced in Proposition
\ref{techest}.
By choosing $\vep$ even smaller if necessary, we can also assume that,
for every $r\in(0,2\vep)$,
%
\begin{equation}
\label{ini2} \bigl|{\bb E}\bigl[F\bigl(\be^{(r)},Z^{(r)}\bigr)
\bigr]- {\bb E}\bigl[F\bigl(\be^{(0)},Z^{(0)}\bigr)\bigr]\bigr| <b.
\end{equation}
In the following, we fix $\vep\in(0,\frac{1}{4})$ so that the
previous two bounds hold.

If $\alpha,\delta\in(0,\frac{1}{4})$, we let $\wh
C^{(\alpha,\delta,n)}$ and $\wh V^{(\alpha,\delta,n)}$ be respectively
the contour and the label function of the labeled tree
$(\t_{(w_n^{\alpha,\delta})},\u_{(w_n^{\alpha,\delta})})$.\vspace*{9pt}

\noindent{\textbf{First step}}. We verify that we can find $\alpha,\delta
\in(0,\vep)$ such that, for all sufficiently large $n$, we have both
$\ov{\bb P}_{n,\mc H}(\wt E_n^{\alpha,\delta})\geq1-b$, and
%
\begin{equation}
\label{concotech0} \bigl| \ov{\bb E}_{n,\mc H}\bigl[ F\bigl(\mc
C^n, \mc V^n\bigr)\bigr] - \ov{\bb E}_{n,\mc H}
\bigl[ {\mathbf1}_{\wt E_n^{\alpha,\delta}} F\bigl(\wh{\mc C}^n,\wh{\mc
V}^n\bigr)\bigr]\bigr| \leq b,
\end{equation}
where similarly as in (\ref{scalingcontour}), we have set, for every
$t\in[0,1]$,
\[
\wh{\mc C}^n_t:= \frac{\sigma}{2}
\bigl(m_{n}^{\alpha,\delta}\bigr)^{-1/2} \wh
C^{(\alpha,\delta,n)}_{2m_{n}^{\alpha,\delta}t},\quad \wh{\mc V}^n_t=
\frac{\sqrt{3}}{2} \sigma^{1/2}\bigl(m_{n}^{\alpha,\delta}
\bigr)^{-1/4} \wh V^{(\alpha,\delta,n)}_{2m_{n}^{\alpha,\delta}t}.
\]

To this end, we use Propositions~\ref{techest} and~\ref{techest2} to
choose $\alpha,\delta\in(0,\vep)$ such that, for all sufficiently
large $n$,
\[
\ov{\bb P}_{n,\mc H}\bigl(\Gamma_n^{\alpha,\delta,\vep}\cap\wt
E_n^{\alpha,\delta}\bigr)\geq1-b/4.
\]
We consider $n$ such that this bound holds and argue on the event
$\Gamma^{\alpha,\delta,\vep}_n\subset E_n^{\alpha,\delta}$. On
this event, the
first visit of the vertex $w^{\alpha,\delta}_n$ by the contour
exploration occurs before time
$2\delta n$ and the last visit of this vertex occurs after time
$2(1-\delta)n$.
From the definition of the pair $(\wh{\mc C}^n,\wh{\mc V}^n)$ we can
find two (random) times $\theta_1\in[0,\delta]$
and $\theta_2\in[1-\delta,1]$, such that
%
\begin{equation}
\label{rescacont} \wh{\mc C}^n_t = \frac{\mc C^n_{\theta_1 +(\theta
_2-\theta_1)t}
-\mc C^n_{\theta_1}}{(\theta_2-\theta_1)^{1/2}},\quad
\wh{\mc V}^n_t = \frac{\mc V^n_{\theta_1 +(\theta_2-\theta
_1)t} -\mc V^n_{\theta_1}}{(\theta_2-\theta_1)^{1/4}}
\quad\forall t \in[0,1].
\end{equation}
It easily follows that
\[
\sup_{t\in[0,\delta]\cup[1-\delta,1]} \bigl|{\mc C}^n_t - \wh{\mc
C}^n_t\bigr| \leq\bigl(1+2(1-2\delta)^{-1/2}\bigr)
\sup_{t\in[0,2\delta]\cup
[1-2\delta,1]} \mc C^n_t \leq4\vep
\]
using the definition of $\Gamma^{\alpha,\delta,\vep}_n$ in the last
inequality.
Still using (\ref{rescacont}), we have also, for every $t\in[0,1]$,
\begin{eqnarray*}
\bigl|{\mc C}^n_{\theta_1 +(\theta_2-\theta_1)t} -\wh{\mc C}^n_{\theta
_1 +(\theta_2-\theta_1)t}\bigr|
&=& \bigl|(\theta_2-\theta_1)^{1/2} \wh{\mc
C}^n_t + \mc C^n_{\theta
_1} - \wh{\mc
C}^n_{\theta_1 +(\theta_2-\theta_1)t}\bigr|
\\
&\leq& \bigl|\wh{\mc C}^n_t -\wh{\mc C}^n_{\theta_1 +(\theta
_2-\theta_1)t}\bigr|
+\bigl(1-(\theta_2-\theta_1)^{1/2}\bigr)\bigl|\wh{\mc
C}^n_t \bigr| + \bigl|\mc C^n_{\theta_1}\bigr|
\end{eqnarray*}
and it follows that
\[
\sup_{t\in[0,1]} \bigl|{\mc C}^n_{\theta_1 +(\theta_2-\theta_1)t} -\wh{\mc
C}^n_{\theta_1 +(\theta_2-\theta_1)t}\bigr| \leq\omega_{\wh{\mc C}^n}(3\delta
)+ 2\delta\bigl\|
\wh{\mc C}^n\bigr\| _\infty+ \vep.
\]
By combining this with the bound on $|\wh{\mc C}^n_t - \mc C^n_t|$
when $t\in[0,\delta]\cup[1-\delta,1]$, we get that
\[
\bigl\|\mc C^n - \wh{\mc C}^n\bigr\|_\infty\leq
\omega_{\wh{\mc C}^n}(3\delta)+ 2\delta\bigl\|\wh{\mc C}^n
\bigr\|_\infty+ 4\vep
\]
on the event $\Gamma^{\alpha,\delta,\vep}_n$. By a similar
argument, we have also
\[
\bigl\|\mc V^n - \wh{\mc V}^n\bigr\|_\infty\leq
\omega_{\wh{\mc V}^n}(3\delta)+ 2\delta\bigl\|\wh{\mc V}^n
\bigr\|_\infty+ 4\vep
\]
on the event $\Gamma^{\alpha,\delta,\vep}_n$.

Now recall that $\ov{\bb P}_{n,\mc H}(\Gamma_n^{\alpha,\delta,\vep
}\cap\wt E_n^{\alpha,\delta})\geq1-b/4$. Since $0\leq F\leq1$, it
follows that
%
\begin{eqnarray}
\label{convcotech11} &&\bigl| \ov{\bb E}_{n,\mc H}\bigl[ F\bigl(\mc
C^n, \mc V^n\bigr)\bigr] - \ov{\bb E}_{n,\mc
H}
\bigl[ {\mathbf1}_{\wt E_n^{\alpha,\delta}} F\bigl(\wh{\mc C}^n,\wh{\mc
V}^n\bigr)\bigr]\bigr|
\nonumber\\[-8pt]\\[-8pt]
&&\quad\leq\frac{b}{2} + \ov{\bb E}_{n,\mc H} \bigl[ \bigl|F\bigl(\mc
C^n,\mc V^n\bigr) - F\bigl(\wh{\mc C}^n,\wh{ \mc V}^n\bigr)\bigr|
{\mathbf1}_{\Gamma_n^{\alpha,\delta,\vep}\cap\wt E_n^{\alpha,\delta}}
\bigr].\nonumber
\end{eqnarray}
From the Lipschitz assumption on $F$ and the preceding bounds on $\|\mc
C^n - \wh{\mc C}^n\|_\infty$ and $\|\mc V^n - \wh{\mc V}^n\|_\infty$,
we see that the second term in the right-hand side is bounded above by
%
\begin{equation}
\label{convcotech1} \ov{\bb E}_{n,\mc H}\bigl[ \bigl(G_\vep
\bigl(\wh{\mc C}^n\bigr)+ G_\vep\bigl(\wh{\mc
V}^n\bigr)\bigr) {\mathbf1}_{\Gamma_n^{\alpha,\delta,\vep}\cap\wt
E_n^{\alpha,\delta}}\bigr].
\end{equation}
The quantity (\ref{convcotech1}) is bounded above by
\[
a_0^{-1} \ov{\bb E}_{n}\bigl[{
\mathbf1}_{E_n^{\alpha,\delta}} \bigl(G_\vep\bigl(\wh{\mc C}^n
\bigr)+ G_\vep\bigl(\wh{\mc V}^n\bigr)\bigr)\bigr]
=a_0^{-1} \ov{\bb E}_{n}\bigl[{
\mathbf1}_{E_n^{\alpha,\delta}} \psi\bigl(n,m_n^{\alpha,\delta}\bigr)\bigr]
\]
where we have used Proposition~\ref{keytech}, and for every integer
$m$ such that
$(1-2\delta)n\leq m\leq n$, we have set
\[
\psi(n,m)= {\bb E}^{\lfloor\alpha\kappa_n\rfloor}_{m} \bigl[G_\vep\bigl({
\mc C}^m\bigr)+ G_\vep\bigl({\mc V}^m\bigr)
\bigr].
\]
We now let $n$ tend to $\infty$.
We note that the ratio $\alpha\kappa_n /\kappa_m$ is bounded above
by $\alpha(1-2\delta)^{-1/4}$ and bounded below by $\alpha$ when $m$
varies over $[(1-2\delta)n,n]\cap\bb Z$. It thus follows from
(\ref{consnakecondi}) that
\[
\limsup_{n\to\infty} \Bigl(\sup_{m\in[(1-2\delta)n,n]\cap\bb Z} \psi(n,m)
\Bigr) \leq\sup_{r\in[\alpha, \alpha(1-2\delta)^{-1/4}]} {\bb E}\bigl
[G_\vep\bigl(
\be^{(r)}\bigr)+G_\vep\bigl(Z^{(r)}\bigr)\bigr] <
\frac{a_0b}{2}
\]
by our choice of $\vep$.
Consequently the quantity (\ref{convcotech1}) is bounded above by
$b/2$ if $n$ is large enough,
and the right-hand side of (\ref{convcotech11}) is then bounded above
by $b$, which gives the bound (\ref{concotech0}).\vspace*{9pt}

\noindent{\textbf{Second step}}. We fix $\alpha$ and $\delta$ as in the
first step above. We then observe that $\wt E^{\alpha,\delta}_n =
E^{\alpha,\delta}_n \cap A_n$, where the event $A_n$ is measurable with
respect to the pair
$(\t^{(w_n^{\alpha,\delta})},\u^{(w_n^{\alpha,\delta})})$. This
measurability property was indeed the motivation for introducing $\wt
E^{\alpha,\delta}_n$. Since the pair $(\wh{\mc C}^n,\wh{\mc V}^n)$ is a
function of $(\t_{(w_n^{\alpha,\delta})},\u_{(w_n^{\alpha,\delta})})$,
and since $\wt E^{\alpha,\delta}_n\subset\mc H$, we can use Proposition~\ref{keytech} to write
\begin{eqnarray*}
\ov{\bb E}_{n,\mc H}\bigl[ {\mathbf1}_{\wt E_n^{\alpha,\delta}} F\bigl(\wh
{\mc
C}^n,\wh{\mc V}^n\bigr)\bigr] &=& \frac{1}{\ov{\bb P}_n(\mc H)} \ov{
\bb E}_{n}\bigl[ {\mathbf1}_{E_n^{\alpha,\delta}} {\mathbf1}_{A_n} F
\bigl(\wh{\mc C}^n,\wh{\mc V}^n\bigr)\bigr]
\\
&=&\frac{1}{\ov{\bb P}_n(\mc H)} \ov{\bb E}_{n}\bigl[ {\mathbf
1}_{E_n^{\alpha,\delta}} {
\mathbf1}_{A_n} \Phi\bigl(n,m^{\alpha,\delta
}_n\bigr)\bigr]
\\
&=&\ov{\bb E}_{n,\mc H}\bigl[ {\mathbf1}_{\wt E_n^{\alpha,\delta}} \Phi
\bigl(n,m^{\alpha,\delta}_n\bigr)\bigr]
\end{eqnarray*}
where, for every integer $m$ such that
$(1-2\delta)n\leq m\leq n$ we have set
\[
\Phi(n,m)=\ov{\bb E}{}^{\lfloor\alpha\kappa_n\rfloor}_{m} \bigl[F\bigl(\mc
C^m,\mc V^m\bigr)\bigr].
\]
If $n$ is large enough, we get from (\ref{consnakecondi}) that
\[
\sup_{(1-2\delta)n\leq m\leq n} \bigl| \Phi(n,m) - {\bb E}\bigl[ F\bigl(\be
^{(\lfloor\alpha\kappa_n\rfloor/\kappa_m)}, Z^{(\lfloor\alpha
\kappa_n\rfloor/\kappa_m)}\bigr)\bigr] \bigr| <b.
\]
Then noting that $\lfloor\alpha\kappa_n\rfloor/\kappa_m \leq2\vep
$ if $(1-2\delta)n\leq m\leq n$, and using (\ref{ini2}), we obtain that
\[
\sup_{(1-2\delta)n\leq m\leq n} \bigl| \Phi(n,m) - {\bb E}\bigl[ F\bigl(\be
^{(0)}, Z^{(0)}\bigr)\bigr] \bigr|<2b,
\]
and we conclude that
\[
\bigl|\ov{\bb E}_{n,\mc H}\bigl[ {\mathbf1}_{\wt E_n^{\alpha,\delta}} F\bigl(\wh
{\mc
C}^n,\wh{\mc V}^n\bigr)\bigr] - {\bb P}_{n,\mc H}
\bigl(\wt E_n^{\alpha,\delta
}\bigr) {\bb E}\bigl[ F\bigl(
\be^{(0)}, Z^{(0)}\bigr)\bigr] \bigr|\leq2b.
\]
By combining this with (\ref{concotech0}), we get
\[
\bigl|\ov{\bb E}_{n,\mc H}\bigl[ F\bigl(\mc C^n, \mc
V^n\bigr)\bigr] - {\bb P}_{n,\mc
H}\bigl(\wt
E_n^{\alpha,\delta}\bigr) {\bb E}\bigl[ F\bigl(\be^{(0)},
Z^{(0)}\bigr)\bigr] \bigr| \leq3b
\]
and finally since $\ov{\bb P}_{n,H}(\wt E_n^{\alpha,\delta})\geq
1-b$, we have
\[
\bigl|\ov{\bb E}_{n,\mc H}\bigl[ F\bigl(\mc C^n, \mc
V^n\bigr)\bigr] - {\bb E}\bigl[ F\bigl(\be^{(0)},
Z^{(0)}\bigr)\bigr] \bigr| \leq4b,
\]
which completes the proof of Theorem~\ref{conv-coding}.

\section*{Acknowledgements}

We are indebted to an anonymous referee for a number of helpful
suggestions. The first author acknowledges financial support from the
Vicerrectorado de investigaci\'on de la PUCP and acknowledges the
hospitality of the D\'epartement de Math\'ematiques d'Orsay where part
of this work was done. The first author's travels and stay were
supported by a grant from the Consejo Nacional de Ciencia, Tecnolog\'ia
e Innovaci\'on Tecnol\'ogica del Per\'u.



\printhistory

\end{document}